\input ./style/arxiv-general.cfg
\documentclass[aos,MSNbibl,seceqn,dvips]{arximspdf}
\makeatletter
   \@ifpackageloaded{graphicx}{}{\usepackage{graphicx}}
\makeatother
\usepackage{mathrsfs}

\doi{10.1214/15-AOS1387}
\volume{44}
\issue{2}
\pubyear{2016}
\firstpage{771}
\lastpage{812}
\docsubty{FLA}

\makeatletter
\newcommand{\rrVert}{\Vert}
\newcommand{\llVert}{\Vert}
\newcommand{\eqref}[1]{(\ref{#1})}
\newtheorem{theorem}{Theorem}
\newtheorem{proposition}{Proposition}
\newtheorem{lemma}{Lemma}
\newproclaim{remark}{Remark}
\newproclaim{assumptions}{Assumptions}

\newcommand{\R}{\mathbb{R}}
\newcommand{\epss}{\varepsilon}
\newcommand{\To}{\Rightarrow}
\newcommand{\tod}{\stackrel{\mathrm{D}}{\to}}
\newcommand{\cov}{\operatorname{cov}}
\makeatother

\begin{document}
\begin{frontmatter}

\title{Amplitude and phase variation of point processes\thanksref{T1}}
\runtitle{Amplitude and phase variation}

\begin{aug}
\author[A]{\fnms{Victor~M.} \snm{Panaretos}\corref{}\ead[label=e1]{victor.panaretos@epfl.ch}}
\and
\author[A]{\fnms{Yoav} \snm{Zemel}\ead[label=e2]{yoav.zemel@epfl.ch}}
\runauthor{V.~M. Panaretos and Y. Zemel}
\thankstext{T1}{Supported by a European Research Council Starting
Grant Award.}
\affiliation{Ecole Polytechnique F\'ed\'erale de Lausanne}
\address[A]{Section de Math\'ematiques\\
Ecole Polytechnique F\'ed\'erale de Lausanne\\
1015 Lausanne\\
Switzerland\\
\printead{e1}\\
\phantom{E-mail:\ }\printead*{e2}}
\end{aug}

%
\received{\smonth{9} \syear{2014}}
%
\revised{\smonth{9} \syear{2015}}

%
\begin{abstract}
We develop a canonical framework for the study of the problem of
registration of multiple point processes subjected to warping, known as
the problem of separation of amplitude and phase variation. The
amplitude variation of a real random function $\{Y(x):x\in[0,1]\}$
corresponds to its random oscillations in the $y$-axis, typically
encapsulated by its (co)variation around a mean level. In contrast, its
phase variation refers to fluctuations in the $x$-axis, often caused by
random time changes. We formalise similar notions for a point process,
and nonparametrically separate them based on realisations of i.i.d.
copies $\{\Pi_i\}$ of the phase-varying point process. A key element in
our approach is to demonstrate that when the classical phase variation
assumptions of Functional Data Analysis (FDA) are applied to the point
process case, they become equivalent to conditions interpretable
through the prism of the theory of optimal transportation of measure.
We demonstrate that these induce a natural Wasserstein geometry
tailored to the warping problem, including a formal notion of bias
expressing over-registration. Within this framework, we construct
nonparametric estimators that tend to avoid over-registration in finite
samples. We show that they consistently estimate the warp maps,
consistently estimate the structural mean, and consistently register
the warped point processes, even in a sparse sampling regime. We also
establish convergence rates, and derive $\sqrt{n}$-consistency and a
central limit theorem in the Cox process case under dense sampling,
showing rate optimality of our structural mean estimator in that case.
\end{abstract}

%

%
\end{frontmatter}

\section{Introduction}\label{sec1}
When analysing the (co)variation of a real random function $\{Y(x):x\in
K\}$ over a continuous compact domain $K$, it can be broadly said that
one may distinguish two layers of variation. The first is \emph
{amplitude variation}. This is the ``classical'' variation that one
would also encounter in multivariate analysis, and refers to the
stochastic fluctuations around a mean level, usually encoded in its
covariance kernel, at least up to second order. In short, this is
variation ``in the $y$-axis.''

The second layer of variation is a non-linear variation peculiar to
continuous domain stochastic processes, and is rarely---if
ever---encountered in multivariate analysis. It arises as the result of random
changes (or deformations) in the time scale (or the spatial domain) of
definition of the process. It can be conceptualised as a composition of
the stochastic process with a random transformation acting on its
domain, or as variation ``in the $x$-axis,'' typically referred to as a
\emph{warp function}. The terminology on amplitude/phase variation is
adapted from trigonometric functions, which may vary in amplitude or phase.

Phase variation arises quite naturally in the study of random phenomena
where there is no absolute notion of time or space, but every
realisation of the phenomenon evolves according to a time-scale that is
intrinsic to the phenomenon itself, and (unfortunately) unobservable.
Processes related to physiological measurements (such as growth curves,
neuronal signals, or brain images), are usual suspects, where phase
variability arises at the level of individual (see the extensive
discussion in Ramsay and Silverman \cite{fdB1,fdB2}); but examples
abound in diverse fields of application of stochastic processes,
perhaps quite prominently in environmental sciences (e.g., Sampson and
Guttorp \cite{guttorp}, and references therein) and pattern recognition
(for instance, handwriting analysis, e.g., Ramsay \cite
{ramsay_handwriting}, or speech analysis, e.g., Hadjipantelis, Aston and Evans
\cite{aston}).

Natural as the confluence of these two types of variation may be,
failing to recognise and correct for their entanglement can obscure or
even entirely distort the findings of a statistical analysis of the
random function (see Section~\ref{sec:fda}). Consequently, it is an
important problem to be able to separate the two, thus correctly
accounting for the distinct contribution of each. If one is able to
only observe a single realisation of the random function $\{Y(x)\}$ in
question, the separation problem is not well-defined unless further
modelling assumptions are introduced. For example, one could assume
that a process should be stationary or otherwise have some invariance
property in the $x$-domain that is measurably perturbed by the phase
variation; and attempt to unwarp it on the basis of this assumption.
Such models can be found in the analysis of random fields (see, e.g.,
Sampson and Guttorp \cite{guttorp}, Anderes and Stein \cite{anderes1},
Anderes and Chatterjee \cite{anderes2}), and of points processes alike
(see, e.g., Schoenberg \cite{schoenberg}, Senoussi, Chadoef and Allard \cite{senoussi}).

In the field of functional data analysis, however, one has the good
fortune of being able to observe multiple i.i.d. realisations $\{
Y_1(x),\ldots,Y_n(x)\}$ of the random function in question. When this is
the case, one may attempt to separate phase and amplitude variation
under less stringent assumptions---in fact in a nonparametric fashion.
Indeed, there is a substantial amount of work on this topic in the
field of functional data, as the problem is in some sense one of the
distinguishing characteristics of FDA as compared to multivariate
statistics (see Section~\ref{sec:fda}).

The purpose of this paper is to investigate the problem of separation
of amplitude and phase variation in the case where one observes
multiple realisations $\{\Pi_1,\ldots,\Pi_n\}$ of \emph{random point
processes} rather than \emph{random functions}. Though the study of
multiple realisations of point processes has been considered prior to
the emergence of FDA (see, e.g., Karr \cite{karr}), treating
realisations of point processes as individual data objects within a
functional data analysis context is a more recent development offering
important advantages; a key paper is that of Wu, M\"uller and Zhang
\cite{wu1} (also see Chiou and M\"uller \cite{chiou} and Chiang, Wang and Huang
\cite{chiang}). Such data may be an object of interest in themselves
(see, e.g., Wu, M\"uller and Zhang \cite{wu1}, Arribas-Gil and M\"uller \cite
{gil}, Wu and Srivastava \cite{wu2}) but may also arise as landmark
data in an otherwise classical functional data analysis (see, e.g.,
Gasser and Kneip \cite{gasser1}, Arribas-Gil and M\"uller \cite{gil}).
The recent surge of interest is exemplified in an upcoming discussion
paper by Wu and Srivastava \cite{EJS}, whose discussion documents early
progress and challenges in the field. One of the main complications
arising in the point process case is that a point processes, when
viewed as a single \emph{datum}, is a discrete random measure. The
nature of such a datum gives rise to different sets of challenges as
compared to FDA. Their ambient space is not a vector space, so point
process variation---whether due to amplitude or phase---is
intrinsically non-linear, calling for an analysis either via a suitable
transformation, or via consideration of an alternative space where
their covariation structure can be suitably analysed. Nevertheless,
this special nature can be seen as a blessing, rather than a curse, as
the case of point processes enjoys important advantages that
considerably simplify the analysis relative to more general functions.

Specifically, we argue that the problem of amplitude and phase
variation in point process data admits a \emph{canonical} framework
through the theory of optimal transportation of measure. Indeed, we
show that this formulation follows unequivocally when employing the
classical phase variation assumptions of functional data analysis to
the point process case (Section~\ref{sec:phase_principles}, Assumptions
\ref{classical_assumptions}). These are proven to be \emph{equivalent}
to a geometrical characterisation of the problem by means of geodesic
variation around a Fr\'echet mean with respect to the Wasserstein
metric (Section~\ref{sec:geometry}, Proposition~\ref
{prop:assumptions}). We show that the special nature of the problem in
the case of point processes renders it identifiable (Section~\ref
{sec:geometry}, Proposition~\ref{prop:uniqueness}) and also allows for
the elucidation of what ``over'' and ``under'' registering means, through
a notion of \emph{unbiased registration} (Section~\ref
{sec:optimality}). We construct easily implementable \emph
{nonparametric} estimators that separate amplitude and phase
(Section~\ref{sec:estimation}) and develop their asymptotic theory,
establishing
consistency in a genuinely nonparametric framework (Section~\ref
{sec:asymptotics}, Theorem~\ref{thm:consistency}) even under sparse
sampling (Remark~\ref{sparse_remark}). In the special case of Cox
processes (randomly warped Poisson processes, see Section~\ref
{sec:cox}), we derive rates of convergence (Theorem~\ref{thm:rate}),
and provide conditions for $\sqrt{n}$-consistency. We also obtain a
central limit theorem for the estimator of the structural mean
(Theorem~\ref{thm:asynorm}), which shows our estimator attains the
optimal rate
under dense sampling and allows for uncertainty quantification
(Remark~\ref{rem:uq}). The finite sample performance methodology is illustrated
by means of examples in Section~\ref{sec:examples}, and a simulation
study in the supplementary material \cite{supplement}.

\section{Amplitude and phase variation of functional data}\label{sec:fda}

In order to motivate our framework for modelling amplitude and phase
variation in point processes, we first revisit the case of functional
data, that is, $n$ independent realisations of a random element of
$L^2[0,1]$, say $\{Y_i(x):x\in[0,1]; i=1,\ldots,n\}$. One typically
understands \emph{amplitude variation} as corresponding to linear
stochastic variability in the observations. That is, assuming that the
mean function is $\mu(x)\in L^2[0,1]$, amplitude variation enters the
model through
\[
Y_i(x)=\mu(x)+Z_i(x),\qquad i=1,\ldots,n,
\]
where the $Z_i(x)$ are mean zero i.i.d. stochastic processes with
covariance kernel $\kappa(s,t)$, typically assumed to be continuous
(equivalently, $Z_i$ are assumed continuous in mean square). In this
setup, the covariation structure of $Y$ can be probed by means of the
Karhunen--Lo\`eve expansion,
%
%
\begin{equation}
\label{uniform_convergence_KL} Y(x)=\mu(x)+\sum_{n=1}^{\infty}
\xi_n\varphi_n(x),
\end{equation}
the optimal Fourier representation of $Y$ in the ortho-normal system of
eigenfunctions of $\kappa$. The equality is understood in $\mathbb
{P}-$mean square, uniformly in $x$. This expansion explains the term
\emph{amplitude variation}: $Y$ varies about $\mu$ by random amplitude
oscillations of the functions $\{\varphi_n\}$. A key feature of this
expansion is the separation of the stochastic component (in the
countable collection $\{\xi_n\}$) and the functional component (in the
deterministic collection $\{\varphi_n\}$).

On the other hand, \emph{phase variation} is understood as the presence
of non-linear variation. Heuristically, this means that there is an
initial random change of time scale, followed by amplitude variation,
yielding \emph{time-warped} curves $\widetilde Y_i$,
%
%
\begin{eqnarray}
\label{gaussian_warping} \widetilde{Y}_i(x)=Y_i
\bigl(T^{-1}_i(x) \bigr)&=&\mu \bigl(T_i^{-1}(x)
\bigr)+Z_i \bigl(T^{-1}_i(x) \bigr)
\nonumber
\\[-8pt]
\\[-8pt]
\nonumber
&=&\mu
\bigl(T_i^{-1}(x) \bigr)+\sum_{n=1}^{\infty}{
\xi_n} {\varphi_n \bigl(T_i^{-1}(x)
\bigr)}.
\end{eqnarray}
The warp functions $T_i:[0,1]\rightarrow[0,1]$ are typically assumed to
be random \emph{increasing functions} independent of the $Z_i$ and with
$\mathbb{E}[T_i(x)]=x$. Consequently, one has
\begin{eqnarray*}
\mathbb{E} \bigl[\widetilde{Y}(x) |T \bigr] &= &\mu \bigl(T^{-1}(x)
\bigr)=\widetilde{\mu}(x);\\
 \operatorname{cov} \bigl\{\widetilde {Y}(x),
\widetilde{Y}(y) \bigr\}
&=& \mathbb{E} \bigl[\kappa \bigl(T^{-1}(x),T^{-1}(y)
\bigr) \bigr] +\operatorname{cov} \bigl\{\widetilde{\mu }(x),\widetilde{\mu}(y)
\bigr\},
\end{eqnarray*}
and thus notices that the right-hand side of equation (\ref
{gaussian_warping}) is no longer interpretable as the Karhunen--Lo\`eve
expansion of $\widetilde{Y}_i$ [the $\varphi_n(T^{-1}(x))$ are \emph
{not} eigenfunctions of the covariance kernel $\operatorname{cov}\{
\widetilde
{Y}(x),\widetilde{Y}(y)\}$]. Indeed, if one ignores phase variation,
and proceeds to analyse the $\widetilde{Y}_i$'s by their own
Karhunen--Lo\`eve expansion, the analysis will be seriously distorted:
the eigenfunctions will be more diffuse and less interpretable (owing
to the effect of attempting to capture horizontal variation via
vertical variation, i.e., local features by global expansions) and the
spectral decay of the covariance operator will be far slower (requiring
the retention of a larger number of components in an eventual principal
component analysis).

The data will then usually come in the form of discrete measurements on
a grid $\{t_j\}_{j=1}^{m}\subset[0,1]$ subject to additive white noise
of variance $\sigma^2>0$,
%
%
\begin{equation}
\label{eq:functional_sampling} \widetilde{y}_{ij}=\widetilde{Y}_i(t_j)+
\varepsilon_{ij}, \qquad i=1,\ldots,n; j=1,\ldots,m,
\end{equation}
assuming of course that the $Y_i$ are continuous. The problem of
separation of amplitude and phase variation can now be seen as that of
recovering the $T_i$ and $Y_i$ from the data $\{\widetilde{y}_{ij}\}
_{i=1}^{n}$, and therefore separating phase variation (fluctuations
of~$T_i$) and amplitude variation (fluctuations of $Y_i$). Doing so
successfully depends on the nature of $T$ (e.g., to guarantee
identifiability), the crystallisation of which is a matter of
assumption. Specifically, more assumptions are needed further to
monotonicity and the expected value being the identity. Indeed, there
does not appear to be a single universally accepted formulation. In
landmark registration, for example, the $T$ are estimated by assuming
that clearly defined landmarks (such as local maxima of the curves or
their derivatives) be optimally aligned across curves (Gasser and Kneip
\cite{gasser1}; see also Gervini and Gasser \cite{selfmodel} for a more
flexible setup). Template methods iteratively register curves to a
template, minimising an overall discrepancy; the template is then
updated, for example, starting from the overall mean (Wang and Gasser
\cite{wang1}; Ramsay and Li \cite{li}). Moment-based registration
proceeds by an alignment of the moments of inertia of the curves (James~\cite{james}). Pairwise separation proceeds by iteratively registering
pairs of observations by means of a penalised sums of square criterion,
and takes advantage of a moment assumption on $T$ being the identity on
average to derive a global alignment (Tang and M\"uller \cite{tang}).
Approaches of a semi-parametric flavour assume a functional form for
$T$ that is known, except for a finite dimensional parameter, and
proceed by likelihood methods in a random-effects type setup (R\o nn
\cite{ronn}; Gervini and Gasser~\cite{gervini+gasser}). Principal
components based registration registers the data so that the resulting
curves have a parsimonious representation by means of a principal
components analysis (the ``least second eigenvalue'' principle; Kneip
and Ramsay~\cite{kneip1}). Elastic registration defines a metric
between curves that is invariant under joint elastic deformation of two
curves by the same warp function, and registers by means of computing
averages with respect to this metric (Tucker, Wu and Sriastava \cite
{fisher-rao}). Multiresolution methods have also been proposed, leading
to the notion of ``warplets'' (Claeskens, Silverman and Slaets \cite{claeskens}).
In recent work, Marron et al. \cite{MBI} consider comparisons between
different registration techniques.

The literature is very rich, and a more in-depth review would be beyond
the scope of the present paper. However, we note that a key conceptual
aspect that recurs in several different estimation approaches in the
literature is the postulate that a registration procedure should
attempt to \emph{minimise phase variability (a fit criterion)} subject
to the constraint that the \emph{registration maps ought to be smooth
and as close to the identity map as possible (a regularity/parsimony
criterion)}. With these key assumptions and principles in mind, we now
turn to consider the case of point process data, and see how these
ideas might be adapted.

\section{Amplitude and phase variation of point processes}

\subsection{Amplitude variation}\label{sec:amp_pp}
Let $\Pi$ be a point process on $[0,1]$, viewed as a random discrete
measure, with the property that
$
\mathbb{E} \{ (\int_0^1d\Pi)^2 \}<\infty$. Defining its mean
measure as
\[
\lambda(A)=\mathbb{E} \bigl\{\Pi(A) \bigr\},\qquad A\in\mathscr{B}
\]
on the collection of Borel sets $\mathscr{B}$ of $[0,1]$, we may
understand amplitude variation as being encoded in the \emph{covariance
measure},
%
%
\begin{equation}
\kappa(A\times B)=\operatorname{cov} \bigl\{\Pi(A),\Pi(B) \bigr\} =\mathbb{E}
\bigl[\Pi(A) \Pi(B) \bigr]-\lambda(A)\lambda(B),
\end{equation}
a signed Radon measure over Borel subsets of $[0,1]^2$. The covariance
measure captures the second order fluctuations of $\Pi(A)$ around its
mean value $\lambda(A)$, as well as their dependence on the
corresponding fluctuations of $\Pi(B)$ around $\lambda(B)$. It
naturally generalises the notion of a covariance operator for
functional data to the case of point process data. Without loss of
generality, we may assume that $\lambda(A)$ is renormalised to be a
probability measure. In the absence of phase variation, estimation of
the covariation structure of $\Pi$ on the basis of $n$ i.i.d.
realisations $\Pi_1,\ldots,\Pi_n$ can be carried out by means of the
empirical versions of $\lambda$ and $\kappa$,
\[
\widehat\lambda_n(A) = \frac{1}{n}\sum
_{i=1}^{n}\Pi_i(A); \qquad \widehat{
\kappa}_n(A\times B) = \frac{1}{n}\sum
_{i=1}^{n}\Pi_i(A)\Pi
_i(B)-\widehat\lambda_n(A)\widehat\lambda_n(B).
\]
These are both strongly consistent (in the sense of weak convergence of
measures with probability 1) as $n\rightarrow\infty$, and in fact one
has the usual central limit theorem in that $\sqrt{n}(\widehat\lambda
_n-\lambda)$ converges in law to a centred Gaussian random measure on
$[0,1]$ with covariance measure $\kappa$ (see, e.g., Karr \cite{karr},
Proposition~4.8).

\subsection{Phase variation: First principles}\label{sec:phase_principles}

Phase variation may be introduced by direct analogy to the functional
case. Assuming that $T_i:[0,1]\rightarrow[0,1]$ are i.i.d. random
homeomorphisms, warped versions of the $\Pi_1,\ldots,\Pi_n$ can be
defined as
\[
\widetilde{\Pi}_i={T_{i}}_{\#}
\Pi_i,\qquad  i=1,\ldots,n,
\]
with ${T_{i}}_{\#}\Pi_i(A)=\Pi_i(T_i^{-1}(A))$ the push-forward of
$\Pi
_i$ through $T_i$. It is natural to assume that the collection $\{T_i\}
$ is independent of the collection $\{\Pi_i\}$. Defining the random
measures $\Lambda_i(A)=\lambda(T^{-1}_i(A))={T_i}_{\#}\lambda(A)$, one
also observes that the conditional mean and covariance measures of $\Pi
_i$ given $T_i$ are
\begin{eqnarray*}
\mathbb{E} \{\widetilde{\Pi} |T \}&=&\Lambda; \\
 \operatorname{cov} \bigl\{
\widetilde{ \Pi}(A),\widetilde{\Pi}(B) \bigr\} &=&\mathbb{E} \bigl\{\kappa
\bigl(T^{-1}(A),T^{-1}(B) \bigr) \bigr\}+\operatorname{cov}
\bigl\{\Lambda(A),\Lambda(B) \bigr\},
\end{eqnarray*}
in analogy to the functional case. Furthermore, if $\Pi
_i([0,t))-\lambda
([0,t))$ is mean-square continuous (equivalently, if $\operatorname
{var}[\Pi
(0,t)]$ is continuous), we have an expansion similar to that of
equation (\ref{uniform_convergence_KL}) for the compensated process,
and the warped compensated process
\begin{eqnarray*}
\Pi_i \bigl([0,t )\bigr)-\lambda \bigl([0,t )\bigr) &=& \sum
_{n=1}^{\infty}\zeta_n \psi_n(t);\\
\widetilde\Pi_i \bigl([0,t )\bigr)-(T_{\#}\lambda)
\bigl([0,t )\bigr) &=& \sum_{n=1}^{\infty}
\zeta_n \psi_n \bigl(T^{-1}(t) \bigr),
\end{eqnarray*}
where $\{\psi_n\}$ are the eigenfunctions of $\kappa(s,t)=\kappa\{
[0,s],[0,t] \}$, in analogy with equation (\ref
{gaussian_warping}). The task of separation of amplitude and phase
variation amounts to constructing estimators $\{\widehat T_i\}$ and $\{
\widehat\Pi_i\}$ of the random maps $T_i$ and of the unwarped
(registered) point processes $\{\Pi_i\}$, respectively, on the basis of
$\widetilde{\Pi}_1,\ldots,\widetilde{\Pi}_n$. Phase variation is then
attributed to the $\{\widehat T_i\}$ and amplitude variation to the $\{
\widehat\Pi_i\}$. As with the case of random curves, if consistent
separation is to be achievable, we will need to impose some basic
assumptions on the precise stochastic and analytic nature of the $\{
T_i\}$. These will come in the form of \emph{unbiasedness} and \emph
{regularity}.

%
\begin{assumptions}\label{classical_assumptions}
The maps $T_i:[0,1]\rightarrow[0,1]$ are i.i.d. random homeomorphisms
distributed as $T$, independently of the point processes $\{\Pi_i\}$.
The random map $T$ satisfies the following two conditions:
\begin{longlist}[(A1)]
\item[(A1)] \textit{Unbiasedness}: $\mathbb{E}[T(x)]=x$ almost
everywhere on $[0,1]$.

\item[(A2)] \textit{Regularity}: $T$ is monotone increasing almost surely.
\end{longlist}
\end{assumptions}

Assumption (A1) asks that the average time change $\mathbb{E}[T(x)]$ be
the identity: on average, the ``objective'' time-scale should be
maintained, so that time is not overall sped up or slowed down. Now,
since $T$ is already a homeomorphism, it is bound to be monotone,
either increasing or decreasing. The regularity assumption (A2) asks
that $T$ represent a proper warping of time (time change): if (A2) were
to fail, we would have a time reversal, which is rather problematic in
most applied settings. Indeed, these assumptions are arguably \emph
{sine qua non} in the classical FDA phase variation literature, perhaps
supplemented with further conditions as discussed earlier. We will now
see that now such further conditions are unnecessary in the point
process case, as they \emph{derive} from the basic assumptions (A1)
and~(A2).

\subsection{Phase variation: Geometry}\label{sec:geometry}

Though our unbiasedness and regularity assumptions stem from first
principles related to warping, they in fact are fully compatible with
an elegant geometrical interpretation of phase variation---indeed one
that opens the way for its consistent separation.

One may consider the space of all diffuse probability measures on
$[0,1]$ as a metric space, endowed with the so-called $L^2$-Wasserstein
distance (also known as Mallows' distance, or earth-mover's distance),
%
%
\begin{equation}
\label{wasserstein_definition} d(\mu,\nu)=\inf_{Q\in\Gamma(\mu
,\nu)}\sqrt{\int
_0^1\bigl|Q(x)-x\bigr|^2\mu(dx)},
\end{equation}
where $\Gamma(\mu,\nu)$ is the collection of mappings
$Q:[0,1]\rightarrow[0,1]$ such that $Q_{\#}\mu=\nu$. The metric $d$ is
related to the so-called Monge problem of optimally transferring the
mass of $\mu$ onto $\nu$, with the cost of transferring a unit of mass
from $x$ to $y$ being equal to their squared distance, $|x-y|^2$. In
the case of diffuse measures $(\mu,\nu)$, the infimum in equation
(\ref
{wasserstein_definition}) is attained at a unique map $T\in\Gamma(\mu
,\nu)$ that is explicitly given by
\[
T=F_\nu^{-1}\circ F_{\mu},
\]
where $F_{\mu}(t)=\int_0^t\mu(dx)$, $F_{\nu}(t)=\int_0^t\nu(dx)$ are
the cumulative distribution functions corresponding to the two
measures, and $F^{-1}_{\nu}$ is the quantile function $F^{-1}_{\nu
}(p)=\inf\{y\in[0,1]:F_\nu(y)\geq p\}$
(see Villani \cite{villani}, Chapter~7;
Bickel and Freedman \cite{freedman}). Consequently, the
optimal map $T$ inherits the regularity properties of the measures $\mu
$ and $\nu$, and does not require any further regularising assumptions.
For example, if both measures admit continuous densities strictly
positive on $[0,1]$, then $T$ is a homeomorphism, but further
smoothness assumptions on the densities will carry over to smoothness
properties of the optimal maps.

When equipped with the metric $d$, the space of all diffuse probability
measures on $[0,1]$ is a \emph{length space} (also known as \emph{inner
metric space}), and the optimal Monge maps $T$, known as optimal
transport maps, generate the geodesic structure of this space.
Specifically, given any diffuse pair $(\mu,\nu)$, there is a unique
geodesic curve $\{\gamma(t):t\in[0,1]\}$ with endpoints $\mu$ and
$\nu$
that is explicitly given by
\[
\gamma(t)= \bigl[t T+(1-t)I \bigr]_{\#}\mu,\qquad t\in[0,1],
\]
where $T$ is the optimal coupling map of $\mu$ and $\nu$, and $I$ is
the identity mapping~\cite{villani}, equation (5.11). The following
proposition demonstrates how this optimal transportation geometry is
inextricably linked with the first principles of phase variation, as
encapsulated in assumptions (A1) and (A2).
%

%
\begin{proposition}\label{prop:assumptions} Let $\lambda$ have strictly
positive density with respect to\break Lebesgue measure on $[0,1]$. A random
map $T:[0,1]\rightarrow[0,1]$ satisfies assumptions \textup{(A1)} and \textup{(A2)}, if
and only if it satisfies assumptions \textup{(B1)} and \textup{(B2)} as stated below:
\begin{longlist}[(B1)]
\item[(B1)] \textit{Unbiasedness}: Given any diffuse probability
measure $\gamma$ on $[0,1]$, we have
\[
\mathbb{E} \bigl\{d^2(T_{\#}\lambda,\lambda) \bigr\}\leq
\mathbb{E} \bigl\{ d^2(T_{\#}\lambda,\gamma) \bigr\}.
\]

\item[(B2)] \textit{Regularity}: Whenever $T_{\#}\lambda=Q_{\#
}\lambda
$, for some homeomorphism $Q:[0,1]\rightarrow[0,1]$, it must be that
\[
\int_0^1 \bigl|T(x)-x \bigr|^2\lambda(dx)
\leq\int_0^1 \bigl|Q(x)-x \bigr|^2
\lambda(dx) \qquad\mbox{almost surely}.
\]
\end{longlist}
\end{proposition}

\begin{figure}

\includegraphics{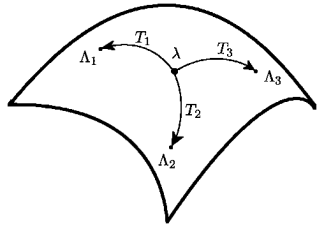}

\caption{Schematic representation of the geometry of phase variation
implied by our assumptions.}
\label{fig:geometry}
\end{figure}

In the optimal transportation geometry, the equivalent assumptions (B1)
and (B2) have a clear-cut interpretation. Assumption (B2) implies that
the conditional means $\Lambda_i={{T}_i}_{\#}\lambda$ of the
warped processes correspond to perturbations of the structural mean
measure $\lambda$ along geodesics (see Figure~\ref{fig:geometry}).
Furthermore, in the presence of (B2), assumption (B1) stipulates that
these geodesic perturbations are ``zero mean'' in that the structural
mean measure $\lambda$ is a Fr\'echet mean of the $\Lambda_i$,
\[
\mathbb{E} \bigl\{d^2(\Lambda,\lambda) \bigr\}\leq\mathbb{E} \bigl\{
d^2(\Lambda,\gamma) \bigr\} \qquad\mbox{for any probability measure }
\gamma.
\]
Notice how these assumptions also mimic the additional estimation
principles encountered in the phase variation of functional data (as
discussed in the end of Section~\ref{sec:fda}): we ask that the warp
maps be such that \emph{phase variability around the structural mean be
minimised \textup{[}our unbiasedness assumption \textup{(B1)]}} subject to the constraint
that the \emph{registration maps deviate as least as possible from the
identity map \textup{[}our regularity assumption \textup{(B2)]}}. In this case, however,
these principles are \emph{equivalent} to the basic assumptions, and do
not have to be added as supplementary.

%

Furthermore, the following proposition establishes that if $\lambda$ is
a Fr\'echet mean of each $\Lambda_i$, then it is the unique such Fr\'
echet mean. Our assumptions, therefore, suffice to guarantee
identifiability of the structural mean (and hence, of the warping
maps). We note that the cumulative distribution function of $\Lambda
=T_{\#}
\lambda$ is strictly increasing almost surely, as a composition of two
such functions.

%
\begin{proposition}[(Identifiability)]\label{prop:uniqueness}
Let $\Lambda$ be a diffuse random probability measure on $[0,1]$ with a
strictly increasing CDF almost surely. Then the minimiser of the functional
\[
\gamma\mapsto\mathbb{E} \bigl\{d^2(\Lambda,\gamma) \bigr\},
\]
defined over probability measures $\gamma$ on $[0,1]$, exists and is unique.
\end{proposition}

\subsection{{Phase variation: Measures vs. densities}}\label
{sec:measures_vs_densities}

One should note that postulating that $\widetilde{\Pi}=T_{\#}\Pi$
induces phase variation of the conditional mean measure relative to the
structural mean measure, $\Lambda=T_{\#}\lambda$. This is \emph{not}
equivalent to phase variation at the level of the \emph{conditional
mean density}, say $f_{\Lambda}$, relative to the \emph{structural mean
density}, say $f_{\lambda}$. Indeed, if $\Lambda=T_{\#}\lambda$ then
\[
f_{\Lambda}(x) = \biggl[\frac{d}{dx} \bigl(T^{-1}(x) \bigr)
\biggr]f_{\lambda} \bigl(T^{-1}(x) \bigr) ,\qquad  x\in[0,1].
\]
Thus, our framework cannot be equivalent to a model that directly
models phase variation at the level of densities, by postulating (say)
that $f_{\Lambda}(x)=f_{\lambda}(T(x))$. In such a model, phase
variation immediately induces further amplitude variation, as the lack
of a correcting factor $\frac{d}{dx} (T^{-1}(x) )$ means that
the new density is no longer a probability density, and thus the total
measure of $[0,1]$ varies as a result of the variation of $T$ (an
overall amplitude variation effect).

An example of phase variation at the level of densities is the model of
Wu and Srivastava \cite{EJS}, where the smoothed point processes are
viewed as random density functions. These are then registered by
employing the (extended) Fisher--Rao metric, using the algorithm of
Srivastava et al. \cite{anuj_arxiv}. The authors of \cite{anuj_arxiv}
argue that the Fisher--Rao approach consistently recovers phase
variation for models of the type
$
f(x) = U \times g(T(x))$,
where $g$ is a deterministic function, $U$ is a real random variable,
and $T$ is the phase map. In the particular case where phase variation
is of densities, the model for the densities becomes
\[
f_{\Lambda}(x) =U \times f_{\lambda} \bigl(T(x) \bigr).
\]
Comparing the last two displayed equations, we see that the two setups
are compatible when the $T$ are assumed to be linear maps. In this
case, unless $T(x)=x$ almost surely, our two conditions (A1) and (A2)
cannot be consolidated: if we require $\mathbb{E}[T(x)]=x$, for a
non-trivial random map (i.e., $\mathbb{P}[\|T-\mathrm{id}\|_{L^2}>0]>0$), then
$T$ cannot be an almost surely strictly increasing homeomorphism on the
finite interval $[0,1]$.

Whether phase variation is formalised at the level of measure or
density is to some extent a modelling decision. However, it is worth
pointing out that if we wish to understand phase variation as the
result of a \emph{non-linear deformation of the underlying space}
(e.g., a smooth deformation of the coordinate system), then the model
postulating $\Lambda=T_{\#}\lambda$ appears to be the natural choice.

\subsection{Phase variation: The (warped) Poisson process case}\label{sec:cox}

Just as Gaussian processes are the archetypal ones in the analysis of
functional data, Poisson processes are so when it comes to point
processes. It is hence worth to briefly consider the effect of phase
variation as encoded in (A1) and (A2) [and their equivalent versions
(B1) and (B2)] on a Poisson process.

Assume that $\Pi$ is a Poisson point process with mean measure
$\lambda
$, and let $\widetilde{\Pi}=T_{\#}\Pi$ be the warped process, as
before. Then, for any disjoint Borel sets $\{A_1,\ldots,A_k\}\subset
\mathscr{B}$, the random variables $\{\widetilde{\Pi}(A_j)\}_{j=1}^k$
are independent conditional on the random warp map $T$. This is because
$\{T^{-1}(A_j)\}_{j=1}^{k}$ must also be disjoint Borel sets, combined
with the fact that $\{\widetilde{\Pi}(A_1),\ldots,\widetilde{\Pi
}(A_k)\}
=\{{\Pi}(T^{-1}(A_1)),\ldots,{\Pi}(T^{-1}(A_k))\}$, with $\Pi$ being
Poisson. Furthermore, for any $A\in\mathscr{B}$,
\[
\mathbb{P} \bigl[\widetilde{\Pi}(A)=k|T \bigr]=\mathbb{P} \bigl[\Pi
\bigl(T^{-1}(A) \bigr)=k|T \bigr]=e^{-\lambda(T^{-1}(A))}\frac{\lambda
^{k}(T^{-1}(A))}{k!}.
\]
In other words, conditional on $T$, the process $\widetilde
{\Pi}(A)$ is Poisson with mean measure $T_{\#}\lambda$. This
establishes that $\widetilde{\Pi}=T_{\#}\Pi$ is distributionally
equivalent to a \emph{Cox process} with directing random measure $T_{\#
}\lambda=\Lambda$. Consequently, our model for phase variation reduces
to asking that the law of the warped point process is that of a Cox
process, where the random directing measure $\Lambda$ is non-linearly
varying with a Fr\'echet mean (with respect to the Wasserstein
distance) equal to the structural mean. Thus, in the Poissonian case,
the compounding of phase and amplitude variation can be viewed as \emph
{double stochasticity}: the phase variation is attributed to the random
directing measure, and the amplitude variation is attributed to the
Poisson fluctuations conditional on the directing measure. It is worth
comparing this with the framework introduced by Wu, M\"uller and Zhang \cite{wu1},
where point processes are modelled as Cox processes whose driving
log-densities are linearly varying functional data.

\section{Estimation} \label{sec:estimation}

\subsection{Overview of the estimation and registration
procedure}\label
{sec:overview}

Armed with the intuition furnished by the geometrical interpretation of
our assumptions, we may now formulate an estimation strategy. Since the
structural mean measure $\lambda$ is the Fr\'echet mean of the random
measures $\Lambda_i={T_i}_{\#}\lambda$ in the Wasserstein metric, the
natural estimator of $\lambda$ would be the \emph{empirical Fr\'
echet--Wasserstein mean} of $\{\Lambda_1,\ldots,\Lambda_n\}$. Of
course, the true $\{\Lambda_i\}$ are unobservable, and instead we
observe the point processes $\{\widetilde{\Pi}_i\}$. However, since
\[
{T_i}_{\#}\lambda= \Lambda_i =\mathbb{E} \{
\widetilde{\Pi}_i |T_i \},
\]
a sensible strategy is to use proxies (estimates) of the $\{\Lambda
_1,\ldots,\Lambda_n\}$ constructed on the basis of $\{\widetilde{\Pi
}_1,\ldots,\widetilde{\Pi}_n\}$, and attempt to use these to approximate
the empirical Fr\'echet--Wasserstein mean. Our procedure will follow
the steps:
\begin{longlist}[1.]
\item[1.] Estimate the random measures $\Lambda_i$. This may be done, for
example, by carrying out classical density estimation on each
$\widetilde{\Pi}_i$, viewed as a point process with mean measure
$\Lambda_i$. Call these estimators $\widehat{\Lambda}_i$, with
corresponding cumulative distribution functions $\widehat{F}_i(t)=\int_0^t\widehat{\Lambda}_i(dx)$.

\item[2.] Estimate $\lambda$ by the empirical Fr\'echet mean of $\widehat
{\Lambda}_1,\ldots,\widehat{\Lambda}_n$ (with respect to the Wasserstein
metric $d$). We call this estimator the \emph{regularised Fr\'
{e}chet--Wasserstein} mean, and denote it by $\widehat{\lambda}$, with
corresponding cumulative distribution function $\widehat{F}(t)=\int_0^t\widehat{\lambda}(dx)$.

\item[3.] Estimate each $T_i$ by the corresponding optimal transportation
map of $\widehat{\lambda}$ onto $\widehat{\Lambda}_i$. In light of the
discussion in the previous section, this is given by $\widehat
{T}_i=\widehat{F}^{-1}_{i}\circ\widehat{F}$. Equivalently, one may
estimate the registration maps by $\widehat{T}^{-1}_i=\widehat
{T}^{-1}_i=\widehat{F}^{-1}\circ\widehat{F}_{i}$.

\item[4.] Register the point processes by pushing them forward through the
registration maps,
%
%
\begin{equation}
\label{eq:registration} \widehat{\Pi}_i={{\widehat{T}}^{-1}_i}{}_{\#}
\widetilde{\Pi}_i,\qquad i=1,\ldots,n.
\end{equation}
\end{longlist}

Of these steps, the last poses no difficulty once the first three have
been carried out. We consider these in more detail in the following
three subsections.

Before doing so, we comment on how these estimators are modified in the
case where the true mean measure is not a probability measure. In this
case, the true measure, say $\mu$, can always be written as $\mu=c
\lambda$, where $c=\mu([0,1])$ and $\lambda$ is a probability measure.
The parameter $c$ can be easily estimated (consistently) by $\widehat
{c}_n=\frac{1}{n}\sum_{i=1}^{n}\widetilde{\Pi}_i([0,1])$
and the remaining estimators can be constructed by normalising the
$\widehat{\Lambda}_i$ to be probability measures (see, e.g.,
Section~\ref{sec:LambdaEstimation}).

\subsection{Estimation of the conditional mean measures}\label
{sec:LambdaEstimation}

The probability measures $\Lambda_i$ can be estimated by various means;
here we will employ kernel density estimation. For $\sigma>0$, let
$\psi
_{\sigma}(x)=\sigma^{-1}\psi(x/\sigma)$, with $\psi$ a smooth symmetric
probability density function strictly positive throughout the real line
and such that $\int x^2\psi(x)\,dx=1$. Let $\Psi$ be the corresponding
distribution function, $\Psi(t)=\int_{-\infty}^{t}\psi(x)\,dx$.

We consider the following smoothing procedure on a set of points
$x_1,\ldots,x_m$. For $y\in[0,1]$, construct a diffuse probability
measure $\mu_y$ on $[0,1]$ with the strictly positive density
\begin{eqnarray}
\psi_\sigma(x-y)+2b_2\psi_\sigma(x-y)\mathbf1\{x>y
\}+2b_1\psi_\sigma(x-y)\mathbf1\{x<y\}+4b_1b_2,\nonumber\\
\eqntext{x\in[0,1],}
\end{eqnarray}
where $b_1=1-\Psi((1-y)/\sigma)$ and $b_2=\Psi
({-y}/\sigma)$. Indeed, integration gives
\begin{eqnarray*}
\int_0^1\psi_\sigma(x-y)
\,dx&=&1-b_1-b_2;\qquad \int_{y}^1
\psi_\sigma(x-y)\,dx= \frac{1}2-b_1;\\
 \int
_0^{y} \psi_\sigma(x-y)\,dx&=&
\frac{1}2-b_2.
\end{eqnarray*}
The intuition behind this construction is the following. First, we
smooth the Dirac measure $\delta_y$ by the kernel $\psi$ around $y$,
and restrict it to $[0,1]$; this yields a measure with total mass
$1-b_1-b_2$. Then we construct the two one-sided versions of $\psi$
around $y$ with total masses $b_1$ and $b_2$, respectively, and again
restrict them to $[0,1]$. The remaining mass, $4b_1b_2$, is distributed
uniformly across $[0,1]$---it does not really matter what we do with
this mass, and we could have re-distributed it in any diffuse way.
Finally, we construct the estimator
%
%
\begin{eqnarray}
\label{hatLambda} \widehat\Lambda_i=\frac{1}{m_i}\sum
_{j=1}^{m_i}\mu_{x_j},\qquad m_i=
\widetilde\Pi_i \bigl([0,1] \bigr),
\nonumber
\\[-8pt]
\\[-8pt]
\eqntext{(\widehat\Lambda_i=
\mbox{Lebesgue measure if }m_i=0),}
\end{eqnarray}
where the $\{x_j\}_{j=1}^{m_i}$ are the points corresponding to
$\widetilde\Pi_i$.

Our construction was slightly more complicated than usual in order to:
(1) ensure that $\widehat{\Lambda}_i$ is everywhere positive on
$[0,1]$; and, (2) allow us to suitably bound the Wasserstein distance
between the smoothed measure and the discrete measure $\widetilde\Pi
_i/\widetilde\Pi_i([0,1])$. Both these properties will be instrumental
in our theoretical results. Indeed, regarding (2), we have the following.

%
\begin{lemma}\label{smoothing_bound}
In the notation of the current section, when $\widetilde\Pi_i([0,1])>0$
and $\sigma\le1/4$, we have the bound
%
%
\begin{eqnarray}
\label{eq:wasssmooth} &&d^2 \bigl(\widehat\Lambda_i,\widetilde
\Pi_i\large/\widetilde\Pi_i \bigl([0,1] \bigr) \bigr)
\nonumber
\\[-8pt]
\\[-8pt]
\nonumber
&&\qquad\le3 \sigma^2 + 4\max \bigl(\Psi(-1/\sqrt\sigma),1-\Psi(1/\sqrt\sigma)
\bigr).
\end{eqnarray}
\end{lemma}

\subsection{Estimation of the structural mean measure}\label
{sec:lambdaEstimation}
Given our discussion in Section~\ref{sec:geometry}, it makes sense to
use an $M$-estimation approach in order to construct an estimator for
$\lambda$. Since $\lambda$ arises as a minimum of the population
functional $M(\gamma)=\mathbb{E}[d^2(\Lambda,\gamma)]$, with
$\Lambda
=T_{\#}\lambda$, we would like to define an estimator by minimising the
sample functional
\[
M_n(\gamma)=\frac{1}{n}\sum_{i=1}^{n}d^2(
\Lambda_i,\gamma).
\]
Unfortunately, the $\{\Lambda_i\}$ are unobservable, so that they need
to be replaced by their estimators \eqref{hatLambda}, leading to the
proxy functional
\[
\widehat{M}_n(\gamma)=\frac{1}{n}\sum
_{i=1}^{n}d^2(\widehat{\Lambda
}_i,\gamma).
\]
If this functional has a unique minimum, then this is the sample Fr\'
echet mean of the $\{\widehat{\Lambda}_i\}$. This type of optimisation
problem rarely admits a closed-form solution. Gangbo and \'{S}wi\c{e}ch
\cite{gangbo} have considered this in the form of a multi-coupling
problem, and Agueh and Carlier \cite{agueh} in the barycentric
formulation given above. They provide general results on existence and
uniqueness (not restricted to the 1-dimensional case), and
characterising equations. Remarkably, in the 1-dimensional case, these
yield an explicit solution. This can also be determined directly, using
elementary arguments: by our assumption on $\{T_i\}$ being
homeomorphisms and $\lambda$ being diffuse, we know that the measures
$\{\Lambda_i\}$ are diffuse measures supported on $[0,1]$ with
probability 1. It follows that (see, e.g., Villani \cite{villani},
Theorem~2.18)
\begin{eqnarray*}
\widehat{M}_n(\gamma)&=&\frac{1}{n}\sum
_{i=1}^{n}d^2(\widehat{\Lambda
}_i,\gamma) =\frac{1}{n}\sum_{i=1}^{n}
\int_0^1\bigl|\widehat{F}_i^{-1}(x)-F^{-1}_{\gamma
}(x)\bigr|^2
\,dx\\
&=& \frac{1}{n}\sum_{i=1}^{n}\bigl\|
\widehat{F}_i^{-1}-F^{-1}_{\gamma
}
\bigr\|^2_{L^2},
\end{eqnarray*}
with $\|\cdot\|_{L^2}$ the usual norm on $L^2[0,1]$. Therefore, if
there exists an optimum of
\[
\widehat{L}_n(Q)=\frac{1}{n}\sum_{i=1}^{n}
\bigl\|\widehat{F}_i^{-1}-Q\bigr\|^2_{L^2}
\]
and this optimum is a valid quantile function, it must be that the
probability measure corresponding to this quantile function is an
optimum of $\widehat{M}_n(\gamma)$. Indeed, $\widehat{L}_n$ does admit
a unique minimum $\bar{Q}$ given by the empirical mean of the $\{
\widehat{F}_i^{-1}\}$,
\[
\bar{Q}(x)=\frac{1}{n}\sum_{i=1}^{n}
\widehat{F}_i^{-1}(x).
\]
Furthermore, $\bar{Q}$ is non-decreasing and continuous, since each of
the $\widehat{F}_i^{-1}$ is so. It is therefore a valid quantile
function [clearly $\bar Q(0)=0$ and $\bar Q(1)=1$]. We conclude that
$\widehat{M}_n(\gamma)$ attains a unique minimum at the measure
\[
\widehat{\lambda}(A)=\int_{A}\frac{d}{dx} \Biggl(
\frac{1}{n}\sum_{i=1}^{n}
\widehat{F}_i^{-1} \Biggr)^{-1}(x)\,dx,
\]
that is, the probability measure with cumulative distribution function
$\widehat{F}= (\frac{1}{n}\sum_{i=1}^{n}\widehat{F}_i^{-1} )^{-1}$.

\subsection{Estimation of the registration maps}\label{sec:warp_estimation}
Once the conditional mean measures $\{\Lambda_i\}$ and the structural
mean measure $\lambda$ have been estimated, we automatically get the
estimators for the warp and registration maps, respectively,
%
%
\begin{equation}
\label{warping_estimator} \widehat{T}^{-1}_i= \Biggl(
\frac{1}{n}\sum_{j=1}^{n}\widehat
{F}_j^{-1} \Biggr)\circ\widehat{F}_i \quad\mbox{and}\quad
\widehat{T}_i= \bigl(\widehat{T}^{-1}_i
\bigr)^{-1}.
\end{equation}
Note here that if $T$ is the optimal transportation map of $\mu$ onto
$\nu$, the change of variables formula immediately implies that
$T^{-1}$ is the optimal transportation map of $\nu$ onto $\mu$.

\subsubsection{Regularity of the optimal maps}\label{sec:regularity}

As was foretold in the end of Section~\ref{sec:phase_principles}, the
estimation of the warp/registration maps did not require additional
smoothness constraints (and by means of tuned penalties) on $T$. Since
$\widehat{T}^{-1}_i= (\frac{1}{n}\sum_{j=1}^{n}\widehat
{F}_j^{-1} )\circ\widehat{F}_i$, we immediately note that the
estimated maps will be as regular as the estimators of $\lambda$ and
$\Lambda_i$ are, or equivalently, as smooth as the $\widehat{F}_j$. It
follows that the smoothness of the estimated maps will be directly
inherited from any smoothness constraints we place on the estimated
mean and conditional mean measures, and will not require the addition
of any further smoothness penalties.

\section{Bias and over-registering}\label{sec:optimality}
Note that our geometrical framework essentially induces a loss function
in the estimation problem for the structural mean,
\[
\mathscr{L}(\lambda,\delta)=d^2(\lambda,\delta),
\]
where $\delta=\delta(\Lambda_1,\ldots,\Lambda_n)$ is a candidate
estimator of $\lambda$. Under this loss function, one can consider the
class of \emph{unbiased estimators of the structural mean} (in the
general sense of Lehmann \cite{lehmann}), that is, estimators $\delta
=\delta({\Lambda}_1,\ldots,{\Lambda}_n)$ satisfying
\[
\mathbb{E}_{\lambda}d^2(\lambda,\delta)=\mathbb{E}_{\lambda}
\mathscr{L}(\lambda,\delta) \leq\mathbb{E}_{\lambda}\mathscr {L}(\gamma,
\delta)=\mathbb{E}_{\lambda
}d^2(\gamma,\delta)
\]
for all diffuse measures $\lambda$ and $\gamma$ on $[0,1]$. A \emph
{biased} estimator $\psi=\psi({\Lambda}_1,\ldots,{\Lambda}_n)$
would be
such that for some measure $\gamma$,
\[
\mathbb{E}_{\lambda}d^2(\lambda,\psi) >\mathbb{E}_{\lambda}d^2(
\gamma,\psi).
\]
Thus, using a biased estimator in order to estimate the warp functions,
may (on average) occasionally produce registrations that appear to be
``successful'' in the sense that the residual phase variation is small;
but on the other hand, they would be registering to the wrong reference
measure (a bias-variance tradeoff). It would thus appear that \emph
{unbiasedness} is a reasonable requirement in this setup, protecting us
against overfitting (or ``over-registering,'' to be more precise).

Interestingly, unbiased estimators can be characterised in terms of
their quantile functions; in particular, the empirical Fr\'
{e}chet mean of $\{\Lambda_1,\ldots,\Lambda_n\}$ is unbiased.

%
\begin{proposition}[(Unbiased estimators)]\label
{prop:unbiasedness}\label
{prop:umvue}
Let $\Lambda_1,\ldots,\Lambda_n$ be i.i.d. random probability measures
on $[0,1]$ with positive density with respect to Lebesgue measure. Let
$\lambda$ be their (unique) Fr\'echet mean in the Wasserstein metric. A
random measure $\delta$ is unbiased for $\lambda$ if and only if its
expected quantile function is the quantile function of $\lambda$, that is,
%
%
\begin{equation}
\mathbb EF_\delta^{-1}(x)=F_\lambda^{-1}(x)\label{eq:unbiased}
\end{equation}
for almost any $x$. In particular, the (unique) empirical Fr\'
echet--Wasserstein mean of $\Lambda_1,\ldots,\Lambda_n$ is an unbiased
estimator of $\lambda$.
\end{proposition}

We can thus interpret our regularised Fr\'echet--Wasserstein estimator
$\hat\lambda$ as \emph{approximately unbiased}, since it is a proxy for
the unobservable empirical Fr\'echet--Wasserstein mean.

\section{Asymptotic theory}\label{sec:asymptotics}
We now turn to establishing the consistency of the estimators
constructed in the previous section, and the rate of convergence of the
estimator of the structural mean. In the functional case, as
encapsulated in equation~(\ref{eq:functional_sampling}), one would need
to assume that the number of observed curves, $n$, as well as the
number of sampled observations per curve, $m$, diverge. Similarly, we
will need to construct a framework for asymptotics where the number of
point processes $n$, and the number of points per observed (warped)
point process, $\int_0^1\widetilde{\Pi}(dx)$, diverge. To allow for
this, we shall assume that the processes $\{\Pi_i\}$ are infinitely divisible.

%
\begin{theorem}[(Consistency)]\label{thm:consistency}
Let $\lambda$ be a diffuse probability measure whose support is
$[0,1]$, and let $\{\Pi_1^{(n)},\ldots,\Pi_n^{(n)}\}_{n=1}^\infty$
be a
triangular array of row independent and identically distributed
infinitely divisible point processes with mean measure $\tau_n \lambda
$, with $\tau_n>0$ a scalar. Let $\{T_1,\ldots,T_n\}$ be independent and
identically distributed random homeomorphisms on $[0,1]$,
stochastically independent of $\{\Pi_i^{(n)}\}$, and satisfying
assumptions \textup{(B1)} and \textup{(B2)} relative to $\lambda$. Let $\widetilde{\Pi
}_i^{(n)}={T_i}_{\#}\Pi_i^{(n)}$, and $\Lambda_i={T_i}_{\#}\lambda=
\tau
_n^{-1}\mathbb{E} \{\widetilde{\Pi}_i^{(n)} |T_i \}$. (We
shall suppress the dependency on $n$, but we notice that, by
construction, $\Lambda_i$ does not depend on $n$.) If $\sigma_n\to0$
and $\tau_n/\log n\to\infty$ as $n\uparrow\infty$, then:
\begin{longlist}[1.]
\item[1.] The conditional mean measure estimators of Section~\ref
{sec:LambdaEstimation} (constructed with bandwidth $\sigma=\sigma_n$)
are Wasserstein-consistent,
\[
d (\widehat{\Lambda}_i,\Lambda_i )\stackrel{p} {
\longrightarrow}0\qquad \mbox{as }n\uparrow\infty, \forall i.
\]
\item[2.] The regularised Fr\'{e}chet--Wasserstein estimator of the
structural mean measure (as described in Section~\ref
{sec:lambdaEstimation}) is strongly Wasserstein-consistent,
\[
d(\widehat{\lambda},\lambda)\stackrel{{a.s.}} {\longrightarrow}0 \qquad\mbox{as }n
\uparrow\infty.
\]

\item[3.] The warp functions and registration maps estimators of
Section~\ref{sec:warp_estimation} are uniformly consistent,
\begin{eqnarray}
\sup_{x\in[0,1]} \bigl|\widehat{T}_i(x)-T_i(x)
\bigr|\stackrel{p} {\longrightarrow}0\quad\mbox{and}\quad\sup_{x\in[0,1]} \bigl|
\widehat{T}_i^{-1}(x)-T_i^{-1}(x) \bigr|
\stackrel{p} {\longrightarrow}0\nonumber\\
\eqntext{\mbox{as }n\uparrow\infty, \forall i.}
\end{eqnarray}

\item[4.] The registration procedure in equation \eqref{eq:registration} is
Wasserstein-consistent,
\[
d \biggl(\frac{{\widehat{\Pi}}_i}{\widehat{\Pi}_i([0,1])},\frac
{\Pi_i}{\Pi
_i([0,1])} \biggr)\stackrel{p} {
\longrightarrow}0 \qquad\mbox{as }n\uparrow\infty, \forall i.
\]
\end{longlist}
Under the additional conditions that $\sum_{n=1}^\infty\tau
_n^{-2}<\infty$ and $\mathbb E [\Pi_1^{(1)}([0,1]) ]^4<\infty
$, the convergence in (1), (3) and (4) holds almost surely.
\end{theorem}

%
\begin{remark}\label{sparse_remark}
The assumption that $\tau_n/\log n\to\infty$ is only needed in order to
avoid empty point processes. It requires that the number of observed
processes should not grow too rapidly relative to the mean number of
points observed per process. This condition can be compared to similar
conditions relating the number of discrete observations per curve in
classical FDA. In a sense, it separates the so-called sparse from the
dense sampling regime (see also Wu, M\"uller and Zhang \cite{wu1}) and
shows that even sparse designs lead to consistency. Notice that no
assumption on the precise rate of convergence of $\sigma_n$ to 0 is
required, and in particular its decay is independent of $\tau_n$.
Indeed, $\sigma_n$ can even be random (e.g.,  sample dependent), provided
it converges to zero in probability (see also Remark~\ref{rem:sigmaRates}).
\end{remark}

%
\begin{remark}
Any (cluster) Poisson process is infinitely divisible, so that this
assumption is not overly restrictive, and allows for the phase varying
point process to be of Cox type, as discussed in Section~\ref{sec:cox}
(as a matter of fact, a point process is infinitely divisible if and
only if its finite dimensional distributions are infinitely divisible;
see Daley and Vere-Jones \cite{daley}, Section~10.2, for a detailed
discussion). It allows us to mathematically translate the increasing
expected number of points per process, to a sort of ``i.i.d.'' sampling
framework more similar to the classical FDA one.
\end{remark}

%
\begin{remark}
In conclusion (4), the random quantity $\widehat{\Pi}_i([0,1])=\Pi
_i([0,1])$ is the number of points observed for the \emph{i}th process.
Normalisation by this factor is a technicality ensuring that the
quantities involved are probability measures (or else the Wasserstein
distance would not be well-defined). The actual distance $d (\frac
{{\widehat{\Pi}}_i}{\widehat{\Pi}_i([0,1])},\frac{\Pi_i}{\Pi
_i([0,1])} )$ only depends on the point patterns themselves, and
not on the normalisation.
\end{remark}

In the case of Cox processes, when the processes are Poisson prior to
warping, if we impose a mild constraint on the decay rate of $\sigma
_n$, we can also establish rates of convergence of the estimator
$\widehat\lambda_n$ of the structural mean measure $\lambda$.

%
\begin{theorem}[(Rate of convergence)]\label{thm:rate}
Assume the conditions of Theorem~\ref{thm:consistency}, and suppose in
addition that the processes $\{\Pi_1^{(n)},\ldots,\Pi_n^{(n)}\}
_{n=1}^\infty$ are Poisson. If the kernel $\Psi$ used for the smoothing
has a finite fourth moment $\int_{-\infty}^\infty x^4\,d\Psi
(x)<\infty$,
then $\widehat\lambda_n$ satisfies
\[
d(\widehat\lambda_n,\lambda) \le O_{\mathbb{P}} \biggl(
\frac{1}{\sqrt n} \biggr) + O_{\mathbb{P}} \biggl(\frac{1}{\sqrt
[4]{\tau_n}} \biggr) +
O_{\mathbb{P}} \Biggl( \frac{1}n\sum_{i=1}^n
\sigma_i^{(n)} \Biggr).
\]
Here, $\sigma_i^{(n)}$ is the bandwidth used for constructing
$\widehat
\Lambda_i$, and it is assumed that $\sigma_n=\max_{1\le i\le
n}\sigma
_i^{(n)}\to0$ in probability.
\end{theorem}

%
\begin{remark}
The first term corresponds to the phase variation, the standard $\sqrt
n$ rate resulting from the approximation of a theoretical expectation
by a sample mean. The second term corresponds to the amplitude
variation. The third term corresponds to the bias incurred by the smoothing.
\end{remark}

Theorem~\ref{thm:rate} allows us to conclude that for $\tau_n\ge
O(n^2)$ and $\max_{1\le i\le n}\sigma^{(n)}_i\le O_{\mathbb
{P}}(n^{-1/2})$ we have $\sqrt{n}$-consistency when dealing with Cox
processes, attaining the optimal rate under dense sampling. Indeed,
even more can be said in the dense sampling regime:

%
\begin{theorem}[(Asymptotic normality)]\label{thm:asynorm}
In addition to the conditions of Theorem~\ref{thm:rate}, assume that
$\tau_n/n^2\to\infty$, $\max_{1\le i\le n}\sigma
^{(n)}_i=o_{\mathbb
{P}}(n^{-1/2})$ and that the density of $\lambda$ is bounded below by a
strictly positive constant. Then $\widehat\lambda_n$ is asymptotically
Gaussian, in the sense that
\[
\sqrt n (S_n-\mathrm{id} ) \stackrel{d} {\longrightarrow} Z \qquad\mbox{in
}L^2 \bigl([0,1] \bigr),
\]
where $S_n$ is the optimal transport map from $\lambda$ to $\widehat\lambda_n$, $\mathrm{id}:[0,1]\rightarrow[0,1]$ is the identity map and
$Z$ is a mean-square continuous Gaussian process with covariance kernel
\[
\kappa(x,y) =\operatorname{cov} \bigl\{T(x),T(y) \bigr\},
\]
for $T$ a random warp map distributed as the $\{T_1,\ldots,T_n\}$.
\end{theorem}

%
\begin{remark}[(Uncertainty quantification)]\label{rem:uq}
Since we have uniformly consistent estimators of the maps $\{T_1,\ldots
,T_n\}$, we can construct an empirical estimate of $\operatorname{cov}
\{
T(x),T(y) \}$, which would allow us to carry out uncertainty
quantification on our structural mean estimate (for example in the form
of pointwise confidence intervals of its CDF).
\end{remark}

%
\begin{remark}\label{rem:sigmaRates}
The statements allow the bandwidth $\sigma_i^{(n)}$ to be random. It
follows from Lemma~\ref{lem:havepoints} that the (minimal) number of
points is of the order $O(\tau_n)$. Consequently, if one chooses the
bandwidth by $\sigma_i^{(n)}=\Pi_i^{(n)}([0,1])^{-\alpha}$ for some
$\alpha>0$, then with probability one, $\sigma_n=\max_{1\le i\le
n}\sigma_i^{(n)}\le O(\tau_n^{-\alpha})$. The condition $\sigma
_n=o_{\mathbb{P}}(n^{-1/2})$ then translates to $\tau_n/n^{1/2\alpha
}\to
\infty$, which automatically holds for $\alpha\ge1/4$ due to the
independent assumption that $\tau_n/n^2\to\infty$. Under Rosenblatt's
rule $\alpha=1/5$, one needs the stronger requirement $\tau
_n/n^{5/2}\to
\infty$ for asymptotic normality to hold.\
\end{remark}

\section{Proofs of formal statements}
\mbox{}
\begin{pf*}{Proof of Proposition~\ref{prop:assumptions}}
We begin by showing that conditions (A2) and (B2) are equivalent in
their own right. Then we will show that subject to (B2) being true,
conditions (A1) and (B1) are equivalent. In the language of optimal
transportation, condition (B2) requires that $T$ should be the optimal
transport map between the diffuse measure $\lambda$ and $T_{\#}\lambda
$. By Brenier's theorem (\cite{villani}, Theorem~2.12), it must be that
$T$ is monotone increasing (as the gradient of a convex function on
$[0,1]$), and thus (A2) is implied. Conversely, assume that (A2) holds
true. We know that there is a unique optimal map between $\lambda$ and
$T_{\#}\lambda$ by $\lambda$ being diffuse. By Brenier's theorem, this
map must be monotone increasing, and hence it must be $T$ itself. This
implies (B2).

Consider now condition (B1), which stipulates that given $\gamma$ a
diffuse measure with everywhere positive density $[0,1]$, we have
\[
\mathbb{E} \bigl\{d^2(T_{\#}\lambda,\lambda) \bigr\} \leq
\mathbb{E} \bigl\{d^2(T_{\#}\lambda,\gamma) \bigr\}.
\]
In the presence of (B2), we know that $T$ is an optimal map. It follows
that the left-hand side is
\[
d^2(T_{\#}\lambda,\lambda) =\int\bigl|T(x)-x\bigr|^2\,d
\lambda.
\]
Keeping this in mind, we focus on the right-hand side. Since $\gamma$
is absolutely continuous, it can be written as $Q_{\#}\lambda$, for
some monotone increasing function $Q$, and in fact $Q$ is the optimal
plan between $\lambda$ and $\gamma$ (since any two diffuse measures
have a unique optimal map, which must be monotone increasing). It
follows that
\[
d^2(T_{\#}\lambda,\gamma)=d^2(T_{\#}
\lambda,Q_{\#}\lambda)=\int\bigl|F^{-1}_{T_{\#}\lambda}(x)-F^{-1}_{Q_{\#
}\lambda}(x)\bigr|^2
\,dx.
\]
Now we note that $F_{T_{\#}\lambda}(x)=F_{\lambda}(T^{-1}(x))$, since
$T$ is increasing, and thus $F_{T_{\#}\lambda
}^{-1}(x)=T(F^{-1}_{\lambda
}(x))$; similarly, $Q$ is increasing too, so $F_{Q_{\#}\lambda
}^{-1}(x)=Q(F^{-1}_{\lambda}(x))$. Consequently,
\begin{eqnarray*}
d^2(T_{\#}\lambda,Q_{\#}\lambda)&=&
\int\bigl|F^{-1}_{T_{\#}\lambda
}(x)-F^{-1}_{Q_{\#}\lambda}(x)\bigr|^2
\,dx= \int\bigl|T \bigl(F^{-1}_{\lambda
}(x) \bigr)-Q
\bigl(F^{-1}_{\lambda}(x) \bigr)\bigr|^2\,dx
\\
&=&\int\bigl|T \bigl(F^{-1}_{\lambda}(x) \bigr)-Q
\bigl(F^{-1}_{\lambda}(x) \bigr)\bigr|^2\frac
{f_{\lambda}(F^{-1}_{\lambda}(x))}{f_{\lambda}(F^{-1}_{\lambda
}(x))}
\,dx,
\end{eqnarray*}
where $f_{\lambda}$ is the density of $\lambda$, which we assumed
earlier to be positive everywhere on $[0,1]$. Now we change variables,
setting $y=F^{-1}_{\lambda}(x)$, and observing that $dx=f_{\lambda
}(y)\,dy$, we have
\[
d^2(T_{\#}\lambda,Q_{\#}\lambda) =
\int\bigl|T(y)-Q(y)\bigr|^2f_{\lambda}(y)\,dy =\int\bigl|T(y)-Q(y)\bigr|^2
\lambda(dy).
\]
As a result of our calculations, we see that, in the presence of (B2),
condition (B1) is equivalent to
\[
\mathbb{E}\int\bigl|T(x)-x\bigr|^2\lambda(dx)\leq\mathbb{E}\int
\bigl|T(x)-Q(x)\bigr|^2\lambda(dx)=\int\mathbb{E} \bigl|T(x)-Q(x)\bigr|^2
\lambda(dx),
\]
for all monotone increasing functions $Q$, where the last equality
follows from Tonelli's theorem. The last condition is satisfied if and
only if $\mathbb{E}[T(x)]=x$, $\lambda$-almost everywhere. Thus, when
$\lambda$ has positive density with respect to Lebesgue measure
everywhere on $[0,1]$, we have established that, if (B2) holds, then
(A1) is equivalent to (B1). This completes the proof.
\end{pf*}

\begin{pf*}{Proof of Proposition~\ref{prop:uniqueness}}
Since $\Lambda$ is diffuse and strictly positive, we may re-express the
functional of interest as
\[
M(\gamma)=\mathbb E \bigl[d^2({\Lambda},\gamma) \bigr] =\mathbb E
\biggl[\int_0^1\bigl|{F}_\Lambda^{-1}(x)-F^{-1}_{\gamma
}(x)\bigr|^2
\,dx \biggr] =\mathbb E\bigl\|{F}_\Lambda^{-1}-F^{-1}_{\gamma}
\bigr\|^2_{L^2},
\]
with $\|\cdot\|_{L^2}$ the usual $L^2$ norm. Therefore, if there exists
an optimum of
\[
L(Q)=\mathbb E\bigl\|{F}_\Lambda^{-1}-Q\bigr\|^2_{L^2},\qquad
Q\in L_2 \bigl([0,1] \bigr)
\]
and this optimum is a valid quantile function, it must be that the
probability measure corresponding to this quantile function is an
optimum of $M(\gamma)$. Indeed, $L$ does admit a unique minimum given
by $\Gamma(x)=\mathbb{E}[F_{\Lambda}^{-1}(x)]$, $x\in[0,1]$, which we
claim is a valid quantile function. Note first that $F^{-1}_{\Lambda}$
is, in fact, a proper inverse of the continuous, strictly increasing
mapping $F_{\Lambda}(x)=\Lambda([0,x])$.
\begin{longlist}[1.]
\item[1.] Since $F^{-1}_{\Lambda}(0)=0$ and $F^{-1}_{\Lambda}(1)=1$ almost
surely, we have $\Gamma(0)=0$ and $\Gamma(1)=1$.

\item[2.] If $x\leq y$, then $F^{-1}_{\Lambda}(x)\leq F^{-1}_{\Lambda}(y)$
almost surely. Consequently,\break $\mathbb{E}[F_{\Lambda}^{-1}(x)]\le
\mathbb
{E}[F_{\Lambda}^{-1}(y)]$ also, proving that $\Gamma$ is non-decreasing.

\item[3.] If $x_k\rightarrow x$ in $[0,1]$, then $X_k=F_{\Lambda
}^{-1}(x_k)\rightarrow F_{\Lambda}^{-1}(x)=X$ almost surely. Since
$|X_k|$ is bounded by 1, the bounded convergence theorem implies that
$\mathbb{E}[X_k]\rightarrow\mathbb{E}[X]$, proving that $\Gamma(x)$ is
continuous at $x$ (and hence everywhere in $[0,1]$ by arbitrary choice
of $x$).\quad\qed
\end{longlist}
\noqed\end{pf*}

\begin{pf*}{Proof of Proposition~\ref{prop:umvue}}
Requiring an estimator $\psi$ to be unbiased translates to
\[
\mathbb{E}_{\lambda}\bigl\|F^{-1}_{\lambda}-F^{-1}_{\psi}
\bigr\|^2_{L^2} \leq\mathbb{E}_{\lambda}
\bigl\|F^{-1}_{\gamma}-F^{-1}_{\psi}
\bigr\|^2_{L^2}.
\]
Since $L^2$ is a linear space, and using Tonelli's theorem to exchange
expectation and integration, the unbiasedness condition is equivalent
to requiring that
\[
\mathbb{E}_{\lambda} \bigl[F^{-1}_{\psi}(x) \bigr]
=F^{-1}_{\lambda}(x) \qquad\mbox{almost everywhere}.
\]
To show that this is indeed the case for the empirical Wasserstein mean
$\delta$, we note that
\[
F^{-1}_{\Lambda_i}=F^{-1}_{(T_i)_{\#}\lambda} =
\bigl(F_{\lambda}\circ T_i^{-1} \bigr)^{-1}=T_i
\circ F^{-1}_{\lambda},
\]
and so, by Proposition~\ref{prop:assumptions}, it follows that
\[
\mathbb{E}_{\lambda} \bigl[F^{-1}_{\Lambda_i}(x) \bigr] =
\mathbb{E}_{\lambda} \bigl[T_i \bigl(F^{-1}_{\lambda}(x)
\bigr) \bigr] =F^{-1}_{\lambda}(x), \qquad i=1,\ldots,n
\]
almost everywhere on $[0,1]$. Since $F^{-1}_{\delta}(x)=n^{-1}\sum
F_{\Lambda_i}^{-1}(x)$ (see Section~\ref{sec:lambdaEstimation}),
$\mathbb{E}_{\lambda} [F^{-1}_{\delta}(x) ]=F^{-1}_{\lambda
}(x)$ also holds a.e., and the unbiasedness of $\delta$ has been established.
\end{pf*}

\begin{pf*}{Proof of Lemma~\ref{smoothing_bound}}
The squared Wasserstein distance is bounded by the cost of sending all
the mass in $\mu_{x_i}$ to $x_i$. The squared distance between $\mu
_{y}$ and $\delta_y$ is
\begin{eqnarray*}
&&\int_0^1(x-y)^2
\psi_\sigma(x-y)\,dx + 2b_1\int_{y}^1(x-y)^2
\psi_\sigma(x-y)\,dx \\
&&\quad{}+ 2b_2\int_0^{y}(x-y)^2
\psi_\sigma(x-y)\,dx + 4b_1b_2\int
_0^1(x-y)^2\,dx
\\
&&\qquad\le(1+2b_1+2b_2)\int_0^1(x-y)^2
\psi_\sigma(x-y)\,dx + 4b_1b_2
\\
&&\qquad\le(1+2b_1+2b_2)\int_\R(x-y)^2
\psi_\sigma(x-y)\,dx + 4b_1b_2
\\
&&\qquad\le3\int_\R x^2\psi_\sigma(x)
\,dx+4b_1b_2 =3\sigma^2+4b_1b_2\qquad
(\mbox{since }b_1+b_2\le1).
\end{eqnarray*}
The reason we needed the one-sided kernels in addition to the standard
two-sided one is that either $b_1$ or $b_2$ can be large (e.g., if
$y=0$, then $b_2=1/2$), but they cannot both be large simultaneously.
Indeed, when $y\ge\sqrt\sigma$, we have $b_2\le\Psi(-1/\sqrt
\sigma)$
and when $1-y\ge\sqrt\sigma$, $b_1\le1-\Psi(1/\sqrt\sigma)$. When
$\sigma\le1/4$, at least one of these possibilities holds, and since
$0\le b_i\le1$, this implies that
\[
b_1b_2\le\max \bigl(\Psi({-1/\sqrt\sigma} ),1-\Psi({1/
\sqrt\sigma} ) \bigr).
\]
This bound holds for any $y\in[0,1]$, and the conclusion follows.
\end{pf*}

In order to prove Theorem~\ref{thm:consistency}, we first need to
eliminate the possibility of having empty point processes (this is the
only reason we assume $\tau_n/\log n\to\infty$). To this aim, we will
use a seemingly unrelated technical result for binomial distributions.

%
\begin{lemma}[(Chernoff bound for binomial distributions)]\label{lem:cherbin}
Let $N\sim B(\tau,q)$, then
\[
\mathbb P (N\le\tau q/2 ) \le\beta^\tau, \qquad \beta=\beta(q) =2
\bigl({(1-q)/(2-q)} \bigr)^{1-q/2}<1.
\]
\end{lemma}

\begin{pf} For any $t\ge0$, we have
\begin{eqnarray*}
\mathbb P \biggl(N\le\frac{\tau q}2 \biggr) &=&\mathbb P \biggl(\exp (-Nt)\ge
\exp \biggl(-t\frac{\tau q}2 \biggr) \biggr) \le\mathbb E\exp (-Nt) \exp
\biggl(t\frac{\tau q}2 \biggr) \\
&= &\biggl[s^{q/2} \biggl(1-q+
\frac{q}s \biggr) \biggr]^\tau,
\end{eqnarray*}
where $s=e^t\ge1$. A straightforward calculation shows that this is
minimised when $s=(2-q)/(1-q)>1$. The objective value at this point,
$\beta$, must be smaller than the objective value at $s=1$, which is 1.
\end{pf}

%
\begin{lemma}[{[Number of points per process is $O(\tau_n)$]}]\label
{lem:havepoints}
If $\tau_n/\log n\to\infty$, then there exists a constant $C_\Pi>0$,
depending only on the distribution of the $\Pi$'s, such that
\[
\liminf_{n\to\infty}\frac{\min_{1\le i\le n}\Pi
_i^{(n)}([0,1])}{\tau_n} \ge C_\Pi\qquad
\mbox{a.s.}
\]
In particular, there are no empty point processes, so the normalisation
is well-defined.
\end{lemma}

\begin{pf}
Let us denote for simplicity by $\Pi_\tau$ ($\tau>0)$ a point process
that follows the same infinitely divisible distribution as $\Pi
_i^{(n)}$, but with mean measure $\tau\lambda$. Let $p$ be the
probability that $\Pi_1$ has no points (clearly, $p<1$, since $\Pi_1$
has one point in average). It follows from the infinite divisibility
that for any rational $\tau$, the probability that $\Pi_\tau$ has no
points is $p^\tau$. By a continuity argument, this can be extended to
any real value of $\tau$: indeed, the Laplace functional of $\Pi_1$
takes the form (Kallenberg \cite{kallenberg}, Chapter~6)
\[
L_1(f)=\mathbb E \bigl[e^{-\Pi_1f} \bigr] =\exp \biggl(-\int
\bigl(1-e^{-\rho f} \bigr)\,d\mu(\rho) \biggr), \qquad f\in F \bigl([0,1] \bigr),
\]
where $F[0,1]$ is the set of Borel measurable functions $f:[0,1]\to
\mathbb R_+$, and $\mu$ is a Radon measure on the set $P([0,1])$. It
follows that $L_\tau(f)$, the Laplace functional of $\Pi_\tau$, is
$(L_1(f))^\tau$ when $\tau$ is rational, which simply corresponds to
multiplying $\mu$ by the scalar $\tau$. By considering the measure
$\tau
\mu$ for any real $\tau$, we obtain $L_\tau(f)=(L_1(f))^\tau$ for any
value of $\tau$. The Laplace functional completely determines the
distribution of the process; in particular, the probability of $\Pi
_\tau
$ having no points is obtained as the limit
\[
\lim_{m\to\infty}L_\tau(m)=\lim_{m\to\infty}
\bigl(L_1(m) \bigr)^\tau=p^\tau,
\]
by the bounded convergence theorem, where $L_\tau(m)=L_\tau(f)$ for the
constant function $f\equiv m$.

Denote the total number of points by $N_i^{(n)}=\Pi_i^{(n)}([0,1])$,
and assume momentarily that the $\tau_n$'s are integers. Then
$N^{(n)}_i$ is the sum of $\tau_n$ i.i.d. integer valued random
variables $X_i$, each having a probability of $p<1$ to equal zero. (In
the Poisson case, $p=e^{-1}$.) Each $X_i$ is larger than $\mathbf1\{
X_i\ge1\}$, which follows a Bernoulli distribution with parameter
$q=1-p$, and $N_i^{(n)}=\sum X_i\ge\sum\mathbf1\{X_i\ge1\}$. It
follows that for any $m$,
\[
\mathbb P \bigl(N_i^{(n)}\le m \bigr) \le\mathbb P
\bigl(B( \tau_n,q)\le m \bigr).
\]
Since $N_i^{(n)}$ are i.i.d. across $i$, specifying $m=\tau_nq/2$ and
applying Lemma~\ref{lem:cherbin} yields
\begin{eqnarray*}
\mathbb P \biggl(\min_{1\le i\le n}N_i^{(n)} \le
\frac{\tau_nq}2 \biggr) &=& 1 - \biggl[1 - \mathbb P \biggl(N_1^{(n)}
\le\frac{\tau_nq}2 \biggr) \biggr]^n \le1- \bigl(1-
\beta^{\tau_n} \bigr)^n \\
&\le&1- \bigl(1-n\beta^{\tau_n}
\bigr),
\end{eqnarray*}
by the Bernoulli inequality $(1-x)^n\ge1-nx$ (valid for $x\le1$ and
$n$ integer; easily proved by induction on $n$). The right-hand side is
$n^{a+1}$ for $a=(\log\beta)\tau_n/\log n$. Since $\tau_n/\log n\to
\infty$ and $\beta<1$, we have $a\to-\infty$ as $n\to\infty$ so
this is
smaller than $n^{-2}$ for sufficiently large $n$. By the
Borel--Cantelli lemma, the result holds for $C_\Pi=q/2$.

If $\tau_n$ is not an integer, then $N_i^{(n)}$ is the sum of $\lfloor
\tau_n\rfloor$ (the largest integer $\le\tau_n$) i.i.d. random
variables $X_i$ with probability $p'=p^{\tau_n/\lfloor\tau_n\rfloor
}\le
p$ to equal zero. Letting $q'=1-p'\ge q$ and observing that $\mathbb
P(B(k,q')\le m)\le\mathbb P(B(k,q)\le m)$ for any $k$ and any $m$ [or
that $\beta(q')\le\beta(q)$], we obtain
\begin{eqnarray}
\mathbb P \biggl(\min_{1\le i\le n}N_i^{(n)} \le
\frac{\lfloor\tau
_n\rfloor q}2 \biggr) \le\mathbb P \biggl(\min_{1\le i\le n}N_i^{(n)}
\le\frac{\lfloor\tau
_n\rfloor q'}2 \biggr) \le n\beta^{\lfloor\tau_n\rfloor} =n^{a+1} ,\nonumber\\
 \eqntext{\displaystyle a=
\log\beta\frac{\lfloor\tau_n\rfloor}{\log n}.}
\end{eqnarray}
We still have $a\to-\infty$ and since $\tau_n/[\tau_n]\to1$, any
$C_\Pi
<q/2$ will qualify. Thus, the lemma holds with $C_\Pi=q/2$.
\end{pf}

%
\begin{remark}
As the proof shows, the condition $\tau_n/\log n\to\infty$ can be
slightly weakened to
\[
\liminf_{n\to\infty}(\tau_n/\log n) > 2/{-\log\beta}
\]
and the lower bound equals 9.75 in the Poisson case.
\end{remark}

\begin{pf*}{Proof of Theorem~\ref{thm:consistency}}
Maintaining the notation $N_i=N_i^{(n)}=\break \Pi_i([0,1])=\widetilde\Pi
_i([0,1])$, we begin by proving (1). Without loss of generality, assume
that $\tau_n$ takes integer values [otherwise, work with $t_n$, the
greatest integer smaller than $\tau_n$, that is, replace $\tau_n$ by
$t_n$ and $\Lambda_i$ by $(\tau_n/t_n)\Lambda_i$]. Let $i$ be a fixed
integer. Since the processes $\{\Pi_i\}$ are infinitely divisible, it
is clear that the $\{\widetilde{\Pi}_i\}$ must be so too. Consequently,
we note that a single realisation of a point process with mean measure
$\tau_n \Lambda_i$ is equivalent in law to a superposition of $\tau_n$
independent and identically distributed processes $\{P_j^{(n)}\}
_{j=1}^{\tau_n}$, each with mean $\Lambda_i$. We can assume that
$P_j^{(n)}$ are constructed as the push-forward through $T_i$ of
independent and identically distributed point processes $Q_j^{(n)}$
with mean measure $\lambda$, that are independent of $T_i$. It follows
that as $n\to\infty$, (e.g., Karr \cite{karr}, Chapter~4)
\[
\frac{1}{\tau_n}\widetilde\Pi_i\stackrel{d} {=}\frac{1}{\tau_n}
\sum_{j=1}^{\tau_n}P_j^{(n)}
\stackrel{w} {\rightarrow} \Lambda_i \qquad\mbox{in probability},
\]
with ``$\stackrel{w}{\rightarrow}$'' denoting weak convergence of
measures. Since $N_i/\tau_n\stackrel{p}{\rightarrow}1$, it follows by
Slutsky's theorem that
%
%
\begin{equation}
\label{eq:LLNforPi} {\widetilde{\Pi}_i\large/N_i}
\stackrel{w} {\rightarrow} \Lambda_i\qquad  \mbox{in probability}.
\end{equation}
As $[0,1]$ is compact, we conclude that this last convergence also
holds in Wasserstein distance \cite{villani}, Theorem~7.12, in
probability. Noting that by \eqref{eq:wasssmooth} and since $\sigma
_n\to
0$ as $n\to\infty$,
\[
\sup_{\Omega} d (\widehat\Lambda_i,{\widetilde
\Pi_i\large/N_i} ) \to0,\qquad  n\to\infty,
\]
an application of the triangle inequality shows that $d(\widehat
\Lambda
_i,\Lambda_i)\stackrel{p}{\rightarrow}0$, establishing claim (1). For
convergence almost surely, we fix $a\in[0,1]$ and set
\[
S_n=\sum_{j=1}^{\tau_n}X_{nj},\qquad
X_{nj}=P_j^{(n)} \bigl([0,a] \bigr)-\Lambda
_i \bigl([0,a] \bigr),\qquad  j=1,\ldots,\tau_n.
\]
One sees that $S_n^4=\varphi(Q_1^{(n)},\ldots,Q_k^{(n)},T_i)$, where
$k=\tau_n$ and
\begin{eqnarray}
\varphi(q_1,\ldots,q_k,f)= \Biggl[ \sum
_{j=1}^kf{}_{\#}q_j \bigl([0,a] \bigr) - f_
{\#} \lambda \bigl([0,a] \bigr) \Biggr] ^ 4,\nonumber\\
 \eqntext{f\in\operatorname{Hom}[0,1];
q_j\in M_R,}
\end{eqnarray}
(where $M_R$ is the collection of Radon measures on $[0,1]$ endowed
with the topology of weak convergence, and $\operatorname{Hom}[0,1]$
is the
space of homeomorphisms of $[0,1]$ endowed with the supremum norm) is a
measurable function (since it is continuous). It is also integrable
because $0\le f_\#\lambda([0,a])\le1$ and $\mathbb E[{T_i}_\#
Q_j^{(n)}([0,a])]^4\le\mathbb E[Q_j^{(n)}([0,1])]^4<\infty$ by the hypothesis.

Since the arguments of $\varphi$ are independent, the proof of \cite{durrett},
Lemma~6.2.1, can be adapted to show that $\mathbb
E[S_n^4|\sigma(T_i)]=g(T_i)$, where (with a slight abuse of notation)
\begin{eqnarray}
g(f)=\mathbb E_Q \bigl[\varphi \bigl(Q_1^{(n)},
\ldots,Q_k^{(n)},f \bigr) \bigr] = \int dq_1\int
dq_2\cdots\int dq_k\varphi(q_1,
\ldots,q_k,f),\nonumber\\
  \eqntext{f\in\operatorname{Hom}[0,1].}
\end{eqnarray}
The same idea shows that for each $j$,
\begin{eqnarray*}
\mathbb E \bigl[X_{nj}|\sigma(T_i) \bigr]&=& \int
dq_j{T_i}_\#q_j \bigl([0,a] \bigr) -
{T_i}_\#\lambda \bigl([0,a] \bigr)\\
&=& \lambda \bigl(T_i^{-1}
\bigl([0,a] \bigr) \bigr) - \lambda \bigl(T_i^{-1}
\bigl([0,a] \bigr) \bigr) = 0.
\end{eqnarray*}
In words, conditional on $\sigma(T_i)$, $\{X_{nj}\}_{j=1}^{\tau_n}$ are
mean zero independent and identically distributed random variables. One
readily verifies that (see the proof of \cite{durrett}, Theorem~2.3.5,
for the details)
\begin{eqnarray*}
\mathbb E \bigl[S_n^4|\sigma(T_i) \bigr] &=&
\sum_{j=1}^{\tau_n}\mathbb E
\bigl[X_{nj}^4| \sigma(T_i) \bigr] + \sum
_{j<l}\mathbb E \bigl[X_{nj}^2X_{nl}^2|
\sigma(T_i) \bigr]\\
& =&\tau_n\mathbb E
\bigl[X_{11}^4| \sigma(T_i) \bigr] + 3
\tau_n(\tau_n-1)\mathbb E \bigl[X_{11}^2X_{12}^2|
\sigma(T_i) \bigr].
\end{eqnarray*}
Taking again expected values and applying Markov's inequality,
\[
\mathbb P \biggl[ \biggl(\frac{S_n}{\tau_n} \biggr)^4>\epss \biggr]
\le\frac{\mathbb E[S_n^4]}{\epss^4\tau_n^4}= \frac{\mathbb\tau_n
\mathbb E[X_{11}^4] + 3\tau_n(\tau_n-1)\mathbb
E[X_{11}^2X_{12}^2]}{\epss^4\tau_n^4}.
\]
The numerator is finite, and the sum over $n$ of the right-hand side
converges when $\sum_n\tau_n^{-2}<\infty$. As $\epss$ is arbitrary,
$S_n/\tau_n\stackrel{\mathrm{a.s.}}\to0$ by the Borel--Cantelli lemma.

Repeating this argument countably many times, we have
\[
\mathbb P \biggl(\frac{\widetilde\Pi_i([0,a])}{\tau_n}-\Lambda _i \bigl([0,a] \bigr)
\to0 \mbox{ for any rational }a \biggr)=1.
\]
If $a$ is irrational, choose $a_k\nearrow a\swarrow b_k$ rational. We
have the inequalities
\begin{eqnarray*}
\frac{\widetilde\Pi_i([0,a])}{\tau_n}-\Lambda_i \bigl([0,a] \bigr) & \le&
\frac
{\widetilde\Pi_i([0,b_k])}{\tau_n}-\Lambda_i \bigl([0,b_k] \bigr)+\Lambda
_i \bigl([0,b_k] \bigr)-\Lambda_i
\bigl([0,a] \bigr);
\\
\frac{\widetilde\Pi_i([0,a])}{\tau_n}-\Lambda_i \bigl([0,a] \bigr) & \ge&
\frac
{\widetilde\Pi_i([0,a_k])}{\tau_n}-\Lambda_i \bigl([0,a_k] \bigr)+\Lambda
_i \bigl([0,a_k] \bigr)-\Lambda_i
\bigl([0,a] \bigr),
\end{eqnarray*}
from which one concludes that almost surely, for any $k$,
\begin{eqnarray*}
-\Lambda_i((a_k,a]) &\le&\liminf_{n\to\infty}
\frac{\widetilde\Pi_i([0,a])}{\tau_n} - \Lambda_i \bigl([0,a] \bigr) \le\limsup
_{n\to\infty}\frac{\widetilde\Pi_i([0,a])}{\tau_n} - \Lambda _i \bigl([0,a]
\bigr)\\
& \le&\Lambda_i\bigl((a,b_k]\bigr).
\end{eqnarray*}
Letting $k\to\infty$, we see that convergence holds for any continuity
point $a$ of $\Lambda_i$. But $\Lambda_i$ is a continuous measure by
construction. One then easily shows the almost sure analogue of \eqref
{eq:LLNforPi} (take $a=1$) and concludes (1) as above.

In order to prove (2), we note that $\lambda$ being a minimiser of the
functional $M(\gamma)=\mathbb{E}[d^2(\Lambda,\gamma)]$ implies that it
must be the unique such minimiser (this follows by Proposition~\ref
{prop:uniqueness}), since $\Lambda=T_{\#}\lambda$ is diffuse and
everywhere positive on $[0,1]$, and $T$ is a homeomorphism. To
establish the purported convergence, we therefore study the convergence
of $\widehat{M}_n(\gamma)=\frac{1}{n}\sum_{i=1}^{n}d^2(\widehat
{\Lambda}
_i,\gamma)$ to $M$, both viewed as being defined over $P([0,1])$, the
space of probability measures supported on $[0,1]$. Using the triangle
inequality, we may interject the functionals
%
%
\begin{equation}
\label{eq:Mn} M_n(\gamma)=\frac{1}{n}\sum
_{i=1}^{n}d^2(\Lambda_i,
\gamma)
\end{equation}
that is, the empirical functional assuming that the $\Lambda_i$ could
be observed; and
%
%
\begin{equation}
\label{eq:Mstar} M_n^*(\gamma)=\frac{1}n\sum
_{i=1}^n d^2 \biggl(\frac{\widetilde{\Pi
}_i}{N_i},
\gamma \biggr)
\end{equation}
(which is well-defined for $n$ sufficiently large by Lemma~\ref
{lem:havepoints}), and write
\[
\bigl|\widehat{M}_n(\gamma)-M(\gamma)\bigr|\le\bigl|\widehat{M}_n(
\gamma)-M_n^*(\gamma)\bigr|+\bigl|M_n^*(\gamma)-M_n(
\gamma)\bigr|+\bigl|M_n(\gamma)-M(\gamma)\bigr|.
\]
We shall show that each of the three terms in the right-hand side
converges to 0 uniformly.

For any three probability measures $\mu$, $\nu$, $\rho$ on $[0,1]$,
one has
%
%
\begin{eqnarray}
\label{eq:defL} d^2(\mu,\nu)&\leq&\sup_{\theta\in P([0,1]^2)}\int
_{[0,1]}\int_{[0,1]}{|x-y|^2} {
\theta(dx\times dy)}
\nonumber
\\[-8pt]
\\[-8pt]
\nonumber
&\le&\sup_{x,y\in
[0,1]}|x-y|^2=1;
\\
\label{eq:unifLip}\bigl|d^2(\mu,\rho)-d^2(\nu,\rho)\bigr|&=&\bigl|d(\mu,
\rho)+d(\nu,\rho)\bigr|\bigl|d(\mu,\rho)-d(\nu,\rho)\bigr|
\nonumber
\\[-8pt]
\\[-8pt]
\nonumber
&\le&2d(\nu,\mu),
\end{eqnarray}
and consequently
\[
\bigl|\widehat M_n(\gamma)-M_n^*(\gamma)\bigr| \le
\frac{1}n\sum_{i=1}^n\biggl \vert
d^2(\widehat \Lambda_i,\gamma)-d^2 \biggl(
\frac
{\widetilde\Pi_i}{N_i},\gamma \biggr) \biggr\vert\le\frac{2}n\sum
_{i=1}^n d \biggl(\widehat\Lambda_i,
\frac{\widetilde\Pi
_i}{N_i} \biggr).
\]
The right-hand side is independent of $\gamma$ and converges to 0 by
application of \eqref{eq:wasssmooth}.

Similarly,
\[
\sup_{\gamma\in P([0,1])}\bigl|M_n(\gamma)-M_n^*(\gamma)\bigr|
\le\frac{2}n\sum_{i=1}^nd \biggl(
\Lambda_i,\frac{\widetilde\Pi
_i}{N_i} \biggr) =\frac{2}n\sum
_{i=1}^nX_{ni} =2\overline
X_n.
\]
Now $X_{ni}$ is a function of $T_i$ and $\Pi_i^{(n)}$, so by
construction they are i.i.d. across $i$. Setting $Y_{ni}=X_{ni}-\mathbb
EX_{ni}$, we obtain mean zero random variables that are i.i.d. across
$i$ and $|Y_{ni}|\le1$ because $0\le X_{ni}\le1$ by \eqref{eq:defL}.
Applying the argument in~\cite{durrett}, Theorem~2.3.5, again, one obtains
\[
\mathbb P \bigl( (\overline X_n-\mathbb E\overline X_n
)^4>\epss \bigr)=\mathbb P \bigl(\overline{Y}_n^4>
\epss \bigr)\le\frac{n\mathbb E
[Y_{ni}^4 ]+3n(n-1)\mathbb E [Y_{ni}^2 ]}{\epss^4n^4}\le\frac
{3}{\epss^4n^2}.
\]
By the Borel--Cantelli lemma and arbitrariness of $\epss>0$, we have
$|\overline X_n-\mathbb E\overline X_n|\stackrel{\mathrm{a.s.}}\to0$. But
$X_{n1}\stackrel p\to0$ as $n\to\infty$ by \eqref{eq:LLNforPi}, and the
bounded convergence theorem yields $\mathbb E[\overline X_n]=\mathbb
E[X_{n1}]\to0$.

Turning to the term $|M_n(\gamma)-M(\gamma)|$, we remark that the
strong law of large numbers yields
\[
M_n(\gamma)\stackrel{\mathrm{a.s.}} {\longrightarrow} M(\gamma),
\]
for all $\gamma$. To upgrade to uniform convergence over $\gamma$,
observe that by \eqref{eq:unifLip}, both $M_n$ and $M$ are 2-Lipschitz.
By compactness of $P([0,1])$, given $\epss>0$, we can choose an $\epss
$-cover $\gamma_1,\ldots,\gamma_k$. For any $\gamma$, we have
$d(\gamma
,\gamma_j)<\epss$ for some $j$, so
\begin{eqnarray*}
\bigl|M_n(\gamma)-M(\gamma)\bigr|&\le&\bigl|M_n(\gamma_j)-M_n(
\gamma)\bigr|+\bigl|M_n(\gamma_j)-M(\gamma_j)\bigr|+\bigl|M(
\gamma_j)-M(\gamma)\bigr|
\\
&\le& 4d(\gamma,\gamma_j)+\bigl|M_n(\gamma_j)-M(
\gamma_j)\bigr|
\\
&\le& 4\epss+\bigl|M_n(\gamma_j)-M(\gamma_j)\bigr|.
\end{eqnarray*}
Taking $n\to\infty$, then $\epss\to0$, we conclude
\[
\sup_{\gamma}\bigl|M_n(\gamma)-M(\gamma)\bigr|\stackrel{\mathrm{a.s.}} {
\longrightarrow}0,\qquad n\to\infty.
\]
Summarising, we have established that $\sup_\gamma|\widehat
M_n(\gamma
)-M(\gamma)|\stackrel{\mathrm{a.s.}}\to0$. Let $\lambda_n$ be a minimiser of
$\widehat M_n$. By compactness of $P([0,1])$, $\lambda_{n_k}\to\mu$,
for some subsequence and some $\mu$. Then $\widehat M_{n_k}(\widehat
\lambda_{n_k})\to M(\mu)$ by the uniform convergence and continuity of
$\widehat M_n$ and $M$. Since $\widehat M_{n_k}(\widehat\lambda
_{n_k})\le\widehat M_{n_k}(\lambda)\to M(\lambda)$, we get $M(\mu
)\le
M(\lambda)$, which, by uniqueness of $\lambda$ as a minimiser of $M$,
implies that $\mu=\lambda$. This establishes $\widehat\lambda
_n\stackrel{\mathrm{a.s.}}\to\lambda$ with respect to the Wasserstein distance.

To prove part (3), let $F$, $G$, $F_n$ and $G_n$ denote the
distribution functions of $\lambda$, $\Lambda_i$, $\widehat\lambda_n$
and $\widehat{\Lambda}_i$, respectively, restricted to $[0,1]$. Since
$F$ and $G$ are continuous functions, we have $F_n\to F$ and $G_n\to G$
pointwise on $[0,1]$ (either in probability or almost surely, depending
on the assumptions). Furthermore, all these functions are strictly
increasing and continuous, thus invertible. Our goal is to show
\[
G_n^{-1}\circ F_n =\widehat T_i
\to T_i =G^{-1}\circ F\qquad \mbox{uniformly on }[0,1].
\]
Lemma~\ref{lem:pointwisetouniform} below shows that it will suffice to
establish pointwise convergence, as uniform convergence will
immediately follow in our current setup. To this aim, we remark that
since $G$ is continuous on a compact set, it maps closed sets to closed
sets. Being a bijection, this implies that $G^{-1}$ is continuous as well.

We proceed by showing that $G_n^{-1}(t)\to G^{-1}(t)$ for $0<t<1$ (this
is obvious when $t\in\{0,1\}$). Let $x$ be the unique number such that
$G(x)=t$ and let $\epss>0$. Then $G_n(x+\epss)\to G(x+\epss)>t$ so that
$x+\epss\ge G_n^{-1}(t)$, at least for $n$ large. Similarly, $x-\epss
\le G_n^{-1}(t)$ for $n$ large and, $\epss$ being arbitrary, we
conclude that $G_n^{-1}(t)\to x=G^{-1}(t)$.

By Lemma~\ref{lem:pointwisetouniform}, $G_n^{-1}$ converges uniformly
to $G^{-1}$ on $[0,1]$, where the latter is (uniformly) continuous.
Given $\epss>0$, let $\delta$ such that $|t-s|\le\delta\To
|G^{-1}(t)-G^{-1}(s)|\le\epss$. When $n$ is large, $\|F_n-F\|_\infty
\le
\delta$ and $\|G_n^{-1}-G^{-1}\|_\infty\le\epss$. Then, for any
$x\in
[0,1]$, $|F_n(x)-F(x)|<\delta$, whence
\begin{eqnarray*}
G_n^{-1} \bigl(F_n(x) \bigr)&\le&
G_n^{-1} \bigl(F(x)+\delta \bigr)\le G^{-1}
\bigl(F(x)+\delta \bigr)+\epss\le G^{-1} \bigl(F(x) \bigr)+2\epss;
\\
G_n^{-1} \bigl(F_n(x) \bigr)&\ge&
G_n^{-1} \bigl(F(x)-\delta \bigr)\ge G^{-1}
\bigl(F(x)-\delta \bigr)-\epss\ge G^{-1} \bigl(F(x) \bigr)-2\epss.
\end{eqnarray*}
In other words, $\|\widehat T_i-T_i\|_\infty\le2\epss$ for any large
enough $n$, and (3) is proven. Since the functions $\widehat T_i$ and
$T_i$ are again strictly increasing, it also follows that $\widehat
T_i^{-1}$ converges to $T_i^{-1}$ uniformly.

Now, we turn to part (4). Recall that
\[
\widehat{\Pi}_i={\widehat{T}^{-1}_i}{}_{\#}
\widetilde{\Pi}_i= \bigl({\widehat{T}^{-1}_i}
\circ{T_i} \bigr){}_{\#}\Pi_i,\qquad i=1,\ldots,n.
\]
It follows that ${\widehat{T}^{-1}_i}\circ{T_i}$ is a transport plan
of $\Pi_i$ onto $\widehat{\Pi}_i$. Consequently,
\[
d^2 \biggl(\frac{{\widehat{\Pi}}_i}{N_i},\frac{\Pi_i}{N_i} \biggr)\leq\int
_0^1\bigl|\widehat{T}^{-1}_i
\bigl(T_i(x) \bigr)-x\bigr|^2\frac{\Pi_i(dx)}{N_i}\leq\sup
_{x\in[0,1]}\bigl|\widehat{T}^{-1}_i
\bigl(T_i(x) \bigr)-x\bigr|^2.
\]
Note, however, that since $T_i\in\operatorname{Hom}[0,1]$,
\begin{eqnarray*}
\sup_{x\in[0,1]}\bigl|\widehat{T}^{-1}_i
\bigl(T_i(x) \bigr)-x\bigr|&=&\sup_{x\in
[0,1]}\bigl|
\widehat{T}^{-1}_i \bigl(T_i
\bigl(T_i^{-1}(x) \bigr) \bigr)-T_i^{-1}(x)\bigr|\\
&=&
\sup_{x\in
[0,1]}\bigl|\widehat{T}^{-1}_i(x)-T^{-1}_i(x)\bigr|,
\end{eqnarray*}
and the latter converges to zero in probability (or almost surely,
depending on the assumptions) as $n\rightarrow\infty$ from part (3).
\end{pf*}

The following elementary result is stated without proof.

%
\begin{lemma}\label{lem:pointwisetouniform}
Let $F_n:[a,b]\to\R$ be non-decreasing and converge pointwise to a
continuous limit function $F$. Then the convergence is uniform.
\end{lemma}

\begin{pf*}{Proof of Theorem~\ref{thm:rate}}
Let $\lambda_n$ be the minimiser of the empirical functional
$M_n(\gamma)
=\frac{1}{n}\sum_{i=1}^{n}d^2(\Lambda_i,\gamma)$. For a probability
measure $\theta\in P([0,1])$, denote its quantile function $F_\theta
^{-1}\in L^2([0,1])$ by $g(\theta)$. Then \cite{villani}, Theorem~2.18,
says that $g$ is an isometry: $d(\theta,\gamma)=\|g(\theta)-g(\gamma)\|
$. Now
\[
\sqrt n \bigl(g(\lambda_n)-g(\lambda) \bigr) =\sqrt n \Biggl(
\frac{1}n\sum_{i=1}^nF_{\Lambda_i}^{-1}-F_\lambda^{-1}
\Biggr).
\]
These are i.i.d. mean zero random elements in $L^2$, whose norm is
bounded by 1. Therefore, the above expression converges in distribution
to a Gaussian limit $\mathit{GP}$ with $\mathbb E\|\mathit{GP}\|^2<\infty$ as $n\to
\infty
$. In particular,
\[
d(\lambda_n,\lambda) =\bigl\|g(\lambda_n)-g(\lambda)\bigr\|
=O_{\mathbb{P}} \bigl(n^{-1/2} \bigr).
\]

The error resulting from approximating $\lambda_n$ by
$\widehat\lambda_n$, the minimiser of $\widehat M_n$, is
\begin{eqnarray*}
\bigl\llVert g(\lambda_n)-g(\widehat\lambda_n) \bigr
\rrVert&=& \Biggl\llVert\frac{1}n\sum_{i=1}^nF_{\Lambda_i}^{-1}-
\frac{1}n\sum_{i=1}^nF_{\widehat\Lambda_i}^{-1}
\Biggr\rrVert\le\frac{1}n\sum_{i=1}^n
\bigl\llVert F_{\Lambda_i}^{-1}-F_{\widehat\Lambda
_i}^{-1}
\bigr\rrVert\\
&=&\frac{1}n\sum_{i=1}^nd(
\Lambda_i,\widehat\Lambda_i),
\end{eqnarray*}
which, by the triangle inequality, is bounded by
\[
\frac{1}n\sum_{i=1}^nd(
\Lambda_i,\widehat\Lambda_i) \le\frac{1}n\sum
_{i=1}^n d \biggl(\Lambda_i,
\frac{\widetilde{\Pi
}_i^{(n)}}{N_i^{(n)}} \biggr)S_{ni} +\frac{1}n\sum
_{i=1}^n d \biggl(\frac{\widetilde{\Pi
}_i^{(n)}}{N_i^{(n)}},\widehat
\Lambda_i \biggr)S_{ni} +\frac{1}n\sum
_{i=1}^nV_{ni},
\]
where $S_{ni}=1-V_{ni}=\mathbf1\{N_i^{(n)}>0\}$. The first term on the
right-hand side corresponds to the amplitude variation, while the
second corresponds to the smoothing bias. The third term was
introduced to accommodate empty processes. The inequality follows from
the convention that $\widehat\Lambda_i$ is Lebesgue measure when
$N_i^{(n)}=0$ and the distance between any two measures is no larger
than one. This term is negligible by Lemma~\ref{lem:havepoints}:
$\mathbb P(\sum V_{ni}=0)\to1$ so this term ``converges'' to 0 at any rate.

Denote the distances of the amplitude variation by $X_{ni}\in[0,1]$.
For fixed $n$, $X_{ni}$ are i.i.d. across $i$. Since
\[
\mathbb P \Biggl(a_n\frac{1}n\sum
_{i=1}^nX_{ni}>\epss \Biggr) \le
\frac{a_n\mathbb E\sum_{i=1}^nX_{ni}}{n\epss} =\frac{a_n\mathbb
EX_{n1}}\epss,
\]
we seek to find the rate at which $\mathbb EX_{n1}$ vanishes. Let $W_1$
denote the 1-Wasserstein distance. Then equations (7.4) and (2.48) in
Villani \cite{villani} and Fubini's theorem imply that
\begin{eqnarray*}
\mathbb EX_{n1}^2& \le&\mathbb ES_{n1}W_1
\biggl(\Lambda_1,\frac{\widetilde{\Pi
}_1^{(n)}}{\widetilde{\Pi}_1^{(n)}([0,1])} \biggr) = \int
_0^1\mathbb E \biggl\vert \Lambda_1
\bigl([0,t] \bigr) - \frac{\widetilde{\Pi
}_1^{(n)}([0,t])}{N_1^{(n)}}\biggr \vert S_{n1}\,dt \\
&= &\int
_0^1\mathbb E\vert B_t\vert\,dt,
\end{eqnarray*}
where $B_t$ is defined by the above equation. Let $t\in[0,1]$ be fixed.
Since $\widetilde\Pi_1^{(n)}$ is a Cox process with random mean measure
$\Lambda_1$, conditional on $\Lambda_1$ and on \mbox{$N_1^{(n)}=k\ge1$}, $B_t$
follows a centred renormalised binomial distribution; $B_t=B(k,q)/k-q$
with $q=\Lambda_1([0,t])$. Since $B_t$ is centred, the conditional
expectation of $B_t^2$ equals its conditional variance, $q(1-q)/k\le
1/(4k)$ (or 0 if $k=0$). This bound is independent of $\Lambda_1$, so
we conclude that $\mathbb EB_t^2|N_1^{(n)}\le\mathbf1\{N_1^{(n)}>0\}
/(4N_1^{(n)})$.

Now $N^{(n)}_1$ follows a Poisson distribution with parameter $\tau_n$.
Note that if $X\sim \operatorname{Poisson}(\theta)$ then $\mathbb EX^{-1}\mathbf1\{
X>0\}\le2/\theta$, which can be seen by applying the inequality
$1/k\le2/(k+1)$ for $k\ge1$:
\[
\sum_{k=1}^\infty\frac{1}ke^{-\theta}
\frac{\theta^k}{k!} \le\sum_{k=1}^\infty2e^{-\theta}
\frac{\theta^k}{(k+1)!} =2\theta^{-1}\sum_{k=1}^\infty
e^{-\theta}\frac{\theta^{k+1}}{(k+1)!} =\frac{2}\theta \bigl(1-e^{-\theta}-
\theta e^{-\theta} \bigr).
\]
Thus, taking expected values again, we conclude that $\mathbb EB_t^2\le
(2\tau_n)^{-1}$ so that the integrand above is $\mathbb E|B_t|\le
(2\tau
_n)^{-1/2}$. It follows that $\mathbb EX_{n1}^2\le(2\tau_n)^{-1/2}$
and so $\mathbb EX_{n1}\le(2\tau_n)^{-1/4}$. Summarising, the
amplitude variation is of order at most $O_\mathbb P(\tau_n^{-1/4})$.

As for the smoothing bias, it has been shown in the proof of
Theorem~\ref{thm:consistency} that each of the summands is bounded by
$G(\sigma_i^{(n)})$, where
\[
G(\sigma) =\sqrt{3\sigma^2 + 4\max \biggl(\Psi \biggl(
\frac{-1}{\sqrt\sigma} \biggr),1-\Psi \biggl(\frac{1}{\sqrt
\sigma} \biggr) \biggr)}.
\]
If (the distribution corresponding to) $\Psi$ has tails of order
$O(t^{-4})$, then the first summand above dominates, so that $G(\sigma
)\le R_\Psi\sigma$ for some finite constant $R_\Psi$ and all $\sigma
\ge
0$, and
\[
\frac{1}n\sum_{i=1}^n d \biggl(
\frac{\widetilde{\Pi
}_i^{(n)}}{N_i^{(n)}},\widehat\Lambda_i \biggr)S_{ni} \le
\frac{1}n\sum_{i=1}^n G \bigl(
\sigma_i^{(n)} \bigr) \le\frac{1}n\sum
_{i=1}^n R_\Psi\sigma_i^{(n)}
=R_\Psi\frac{1}n\sum_{i=1}^n
\sigma_i^{(n)}.
\]
The result now follows from $d(\widehat\lambda_n,\lambda)\le
d(\widehat
\lambda_n,\lambda_n)+d(\lambda_n,\lambda)$.
\end{pf*}

\begin{pf*}{Proof of Theorem~\ref{thm:asynorm}}
The conditions of the theorem imply that $\sqrt n(g(\widehat\lambda
_n)-g(\lambda_n))$ converges weakly to 0, so that
\[
\sqrt n \bigl(F^{-1}_{\widehat\lambda_n}-F^{-1}_\lambda
\bigr) =\sqrt n \bigl(g(\widehat\lambda_n)-g(\lambda) \bigr) \tod \mathit{GP},
\]
where $\mathit{GP}$ is the Gaussian process defined above. So the first
statement follows from Slutsky's theorem. The assumption that the
density of $\lambda$ is positively bounded below implies that
$u=F_\lambda$ satisfies the hypothesis of Lemma~\ref{lem:composition}
stated after the end of the proof, so that right composition is
continuous on $L^2[0,1]$. By the continuous mapping theorem
\[
\sqrt n(S_n-\mathrm{id}) =\sqrt n \bigl(F^{-1}_{\widehat\lambda_n}\circ
F_\lambda-F^{-1}_\lambda\circ F_\lambda \bigr)
= \bigl[\sqrt n \bigl(F^{-1}_{\widehat\lambda_n}-F^{-1}_\lambda
\bigr) \bigr]\circ F_\lambda\tod \mathit{GP}\circ F_\lambda,
\]
where $S_n$ is the optimal map from $\lambda$ to $\widehat\lambda_n$.

Now $Z=\mathit{GP}\circ F_\lambda$ is also the weak limit of the process
\[
\sqrt n \Biggl(\frac{1}n\sum_{i=1}^nF_{\Lambda_i}^{-1}
\circ F_\lambda- F_\lambda^{-1}\circ F_\lambda
\Biggr) =\sqrt n \Biggl(\frac{1}n\sum_{i=1}^nT_i
- \mathrm{id} \Biggr),
\]
where $T_i$ is the random warp function from $\lambda$ to $\Lambda_i$.
Since these are i.i.d. elements in $L^2$, we see that the covariance of
$Z$ is $\mathbb{E}(T-\mathrm{id})\otimes(T-\mathrm{id})$, that is, the kernel is
\[
\kappa(s,t) = \mathbb E \bigl[ \bigl(T(s)-s \bigr) \bigl(T(t)-t \bigr) \bigr] =
\cov \bigl(T(s),T(t) \bigr),\qquad s,t\in[0,1].
\]
It easily follows from $Z(t)=\mathit{GP}(F_\lambda(t))$ that $Z$ is a Gaussian process.
\end{pf*}

%
\begin{lemma}[(Composition and continuity)]\label{lem:composition}
Let $u:[0,1]\to[0,1]$ be strictly increasing piecewise continuously
differentiable. Suppose that the derivative of $u$ is bounded below by
$\delta>0$. Then the composition from the right $f\mapsto f\circ u$
from $L^p[0,1]$ takes values in $L^p[0,1]$ and it is $\delta^{-1/p}$-Lipschitz.
\end{lemma}

\begin{pf}
Since composition from the right is linear, it is sufficient to prove
continuity around zero. This follows from the change of variables formula
\begin{eqnarray*}
\|f\circ u\|^p &=&\int_0^1\bigl|f^p
\bigl(u(s) \bigr)\bigr|\,ds =\int_{u(0)}^{u(1)}\bigl|f^p(t)\bigr|
\frac{1}{u'(u^{-1}(t))}\,dt \le\frac{1}\delta\int_{u(0)}^{u(1)}\bigl|f^p(t)\bigr|
\,dt\\
& \le&\frac{1}\delta\|f\|^p,
\end{eqnarray*}
since $0\le u(0)\le u(1)\le1$. The statement for $p=\infty$ holds
trivially without any assumptions on $u:[0,1]\to[0,1]$.
\end{pf}

\section{Illustrative examples}\label{sec:examples}

In order to illustrate the estimation framework put forth in the
previous sections, we consider two scenarios involving warped Poisson
processes (equivalently, Cox processes, see Section~\ref{sec:cox}).
More detailed simulations, including comparisons with the Fisher--Rao
approach \cite{anuj_arxiv}, may be found in the supplementary material \cite{supplement}.

\subsection{Explicit classes of warp maps}\label{sec:sine_warps}

\begin{figure}

\includegraphics{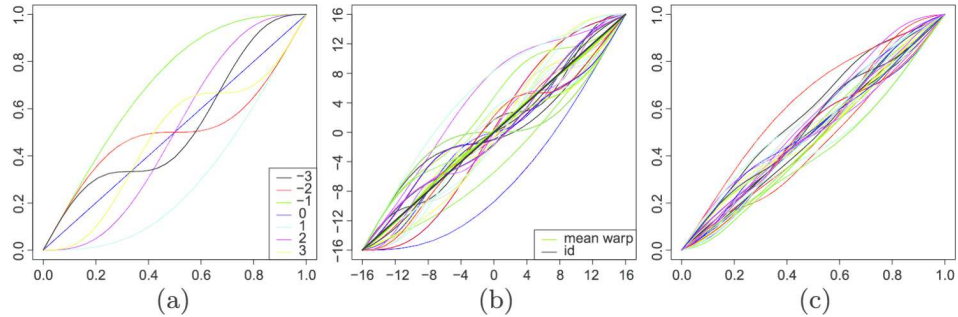}

\caption{\textup{(a)} The functions $\{\zeta_k\}$ for $|k|\le3$;
\textup{(b)}
Realisations of $T$ defined as in equation \protect\eqref
{eq:sinewarps} with
$J=2$ and $K_j\stackrel{d}{=}V_1V_2$ where $V_1$ is Poisson with mean
3, and $\mathbb{P}[V_2=+1]=\mathbb{P}[V_2=-1]=1/2$, independently of
$V_1$; \textup{(c)} Realisations of $T$ defined as in equation
\protect\eqref
{eq:sinewarps} with $J=10$ and $K_j$ as in \textup{(b)}.}\label{fig:sinwarps}
\end{figure}

We first introduce a flexible mixture class of warp maps that provably
satisfies assumptions (A1) and (A2). This can be seen as an extension
of the class considered by Wang and Gasser in \cite{wang1,wang2}. Let
$k$ be an integer and define $\zeta_k:[0,1]\to[0,1]$ by
%
%
\begin{equation}
\label{eq:zetafunction} \zeta_0(x) = x, \qquad\zeta_k(x) = x -
\frac{\sin(\pi kx)}{|k|\pi},\qquad k\in\mathbb{Z}\setminus\{0\}.
\end{equation}
These are strictly increasing smooth functions satisfying $\zeta
_k(0)=0$ and $\zeta_k(1)=1$ for any $k$. Plots of $\zeta_k$ for
$|k|\le
3$ are presented in Figure~\ref{fig:sinwarps}(a). These maps can be made
random by replacing $k$ by an integer-valued random variable $K$. If
the distribution of $K$ is symmetric (around 0), then it is
straightforward to see that
\[
\mathbb{E} \bigl[\zeta_K(x) \bigr] = x \qquad \forall x\in[0,1].
\]
This discrete family of random maps can be made continuous by means of
mixtures: for $J>1$ let $\{K_j\}_{j=1}^{J}$ be i.i.d. integer-valued
symmetric random variables, and $\{U_{(j)}\}_{j=1}^{J-1}$ be the order
statistics of $J-1$ i.i.d. uniform random variables on $[0,1]$,
independent of $\{K_j\}_{j=1}^{J}$. The random map
%
%
\begin{eqnarray}
\label{eq:sinewarps} T(x)=U_{(1)}\zeta_{K_1}(x)+\sum
_{j=2}^{J-1}(U_{(j)}-U_{(j-1)})\zeta
_{K_j}(x)+ (1-U_{(J-1)} )\zeta_{K_J}(x),
\nonumber
\\[-8pt]
\\[-8pt]
\eqntext{x\in[0,1],}
\end{eqnarray}
satisfies assumptions (A1) and (A2). The parameter $J$ can be seen as
controlling the \emph{variance} of $T$: the larger $J$ is, the more
variables are being averaged, and so a law of large numbers effect
yields maps that deviate only slightly from the identity [see
Figure~\ref{fig:sinwarps}(b) and \ref{fig:sinwarps}(c)].

%

%

\subsection{Bimodal Cox processes}
%
%

\begin{figure}

\includegraphics{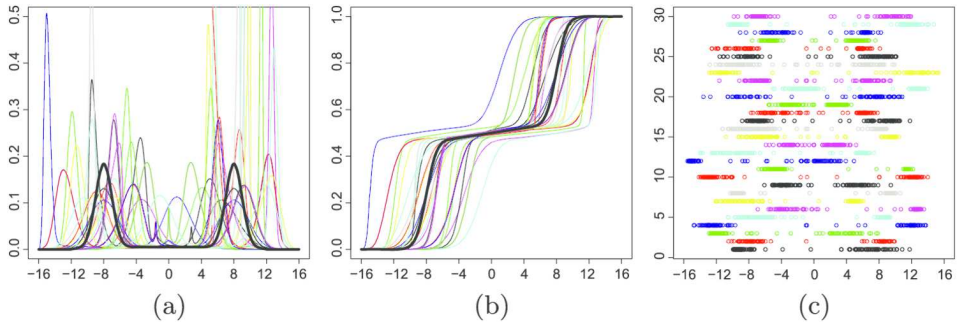}

\caption{\textup{(a)} Thirty warped bimodal densities, with the
structural mean
bimodal density $93\times f$ in solid black; \textup{(b)} Their corresponding
distribution functions, with the structural mean distribution function
in solid black; \textup{(c)} Thirty Cox processes, constructed as
follows: first
we generate $\Pi_i$ as i.i.d. Poisson processes with mean density $f$,
then we warp them by forming $T_{\#}\Pi_i$, where $T$ are the maps
appearing in Figure \protect\ref{fig:sinwarps}\textup{(b)}.}
\label{fig:bidendist}
\end{figure}

We first focus on a scenario where assumptions (B1) and (B2) hold true.
We consider a structural mean measure that is a mixture of three
independent components: two Gaussian distributions (of unit variance),
restricted to the interval $[-16,16]$, and a beta background with
parameters $(1.5,1.5)$, restricted on the interval $[-12,12]$. We wish to
discern the two clear modes (located at $\pm8$), but these may be
smeared by phase variation. The structural mean density is
\[
f(x) =\frac{1-\varepsilon}2 \bigl[\varphi(x - 8) + \varphi(x + 8) \bigr] +
\frac{\varepsilon
}{24} \beta_{1.5,1.5} \biggl(\frac{x+12}{24} \biggr),
\]
where $\varphi$ denotes a standard Gaussian density, $\beta_{\alpha
,\beta}$ is the $\operatorname{Beta}(\alpha,\beta)$ density, and
$\varepsilon
=0.1$ is the strength of the background. We generated 30 independent
Poisson processes with this structural mean measure and $\tau=93$, and
warped them by means of 30 independent warp maps $\{T_i\}$, obtaining
30 warped point processes [Figure~\ref{fig:bidendist}(c)]. The warp maps
$\{T_k\}$ are affinely transformed versions of the maps shown in
Figure~\ref{fig:sinwarps}(b) according to the mapping
\[
g(x)\mapsto32 g \biggl(\frac{x+16}{32} \biggr) - 16
\]
in order to re-scale their support to $[-16,16]$. Recall that the warp
maps in Figure~\ref{fig:sinwarps}(b) were generated using the
definition in equation \eqref{eq:sinewarps}, taking $J=2$ and $K_j$ are
i.i.d., distributed as $V_1V_2$, where $V_1$ is Poisson with mean 3,
and $\mathbb{P}[V_2=+1]=\mathbb{P}[V_2=-1]=1/2$, independently of
$V_1$. These correspond to rather violent phase variation, as can be
seen by the plots of the conditional density/distribution of the warped
processes given the corresponding $T_i$ in Figure~\ref{fig:bidendist}(a)
and \ref{fig:bidendist}(b).

%

\begin{figure}

\includegraphics{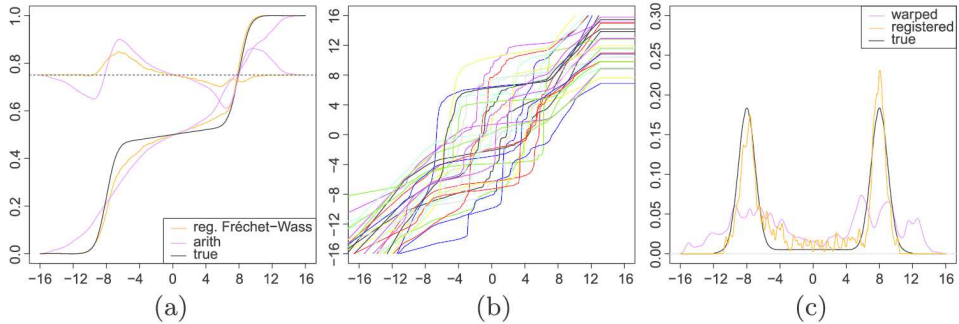}

\caption{\textup{(a)} The empirical arithmetic mean, our estimated regularised
Fr\'echet--Wasserstein mean, and the true mean CDF (the curves
oscillating about the horizontal line $y=3/4$ are residual curves,
centred at $3/4$); \textup{(b)} The estimated warp functions; \textup
{(c)} Kernel estimates
of the density function of the true structural mean, based on the
original spike trains, and on the registered spike trains.}\label{fig4}
\end{figure}

\begin{figure}[b]

\includegraphics{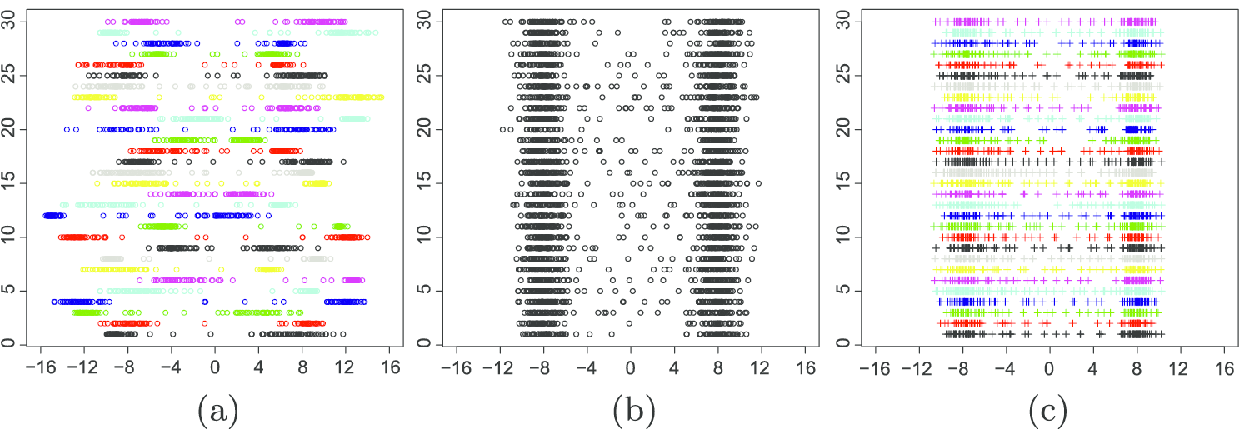}

\caption{Bimodal Cox processes: \textup{(a)} The warped point
processes; \textup{(b)} The
original point processes; \textup{(c)} The registered point processes.}
\label{registration_plots_bimodal}
\end{figure}

\begin{figure}

\includegraphics{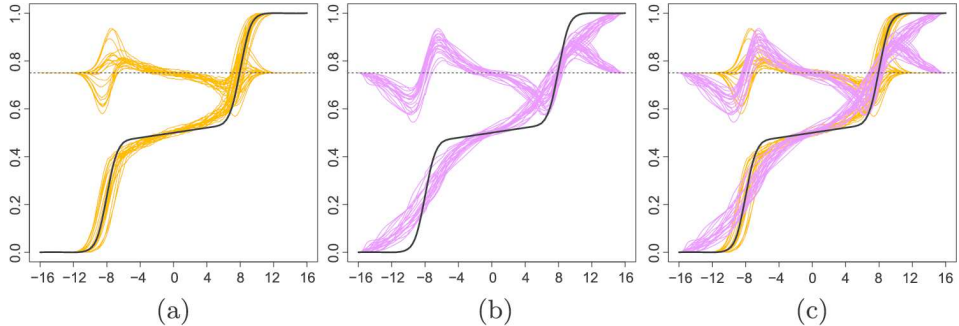}

\caption{\textup{(a)} Comparison of our estimated regularised Fr\'
echet--Wasserstein mean, and the true mean CDF, for 20 independent
replications of the experiment; \textup{(b)} Comparison of the
arithmetic mean,
and the true mean CDF, for the same 20 replications; \textup{(c)} Superposition
of \textup{(a)} and \textup{(b)}. In all three cases, the curves
oscillating about the
horizontal line $y=3/4$ are residual curves, centred at $3/4$.}
\label{bimodal_replications}
\end{figure}

Using the 30 warped spike trains depicted in Figure~\ref{fig:bidendist}(c), we construct the ``regularised Fr\'
echet--Wasserstein'' estimator as described in Section~\ref
{sec:estimation}. A slight deviation is that we use a Gaussian kernel
with bandwidth chosen by unbiased cross validation, rather than the
special kernels developed for the asymptotic theory (with no essential
effect on finite sample performance). We thus obtain estimates of the
warp maps $\{\widehat T_i\}_{i=1}^{30}$ (using the definitions in
Section~\ref{sec:warp_estimation}), depicted in Figure~\ref{fig4}(b), which can be used to register the point processes
(Figure~\ref{registration_plots_bimodal}). The final estimate of the structural
mean distribution function (the regularised Fr\'{e}chet--Wasserstein
estimator) is depicted in Figure~\ref{fig4}(a), and contrasted with
the true structural CDF, as well as with the naive estimate produced by
ignoring warping and averaging the empirical distributions across
trains. We notice that the regularised Fr\'echet--Wasserstein estimator
performs quite well at discerning the two modes of the structural mean
measure, in contrast with the naive estimator which seems to fail to
resolve them. This effect is more clearly portrayed in Figure~\ref
{fig4}(c), which plots kernel estimators of
structural mean density constructed using the original (warped) point
processes, and the registered point processes. It is important to
remark that the minor fluctuations in the density estimate observed are
\emph{not} related to our method of estimation, but are due to the
sampling variation of the spike trains (i.e., they are not intrinsic
to our registration procedure, but to the kernel density estimation
procedure), and could be reduced by more careful choice of bandwidth.
Figure~\ref{bimodal_replications} presents the sampling variation of
the regularised Fr\'echet--Wasserstein estimator, and contrasts it with
the sampling variation of the naive arithmetic estimator for 20
independent replications of the same experiment. We notice that the
naive estimator is clearly biased in the neighbourhoods around the two
peaks, and appears to fluctuate around a straight line. In contrast,
the smoothed Fr\'echet mean---though presenting fluctuations around the
two peaks---appears approximately unbiased. Indeed, its variation is
very clearly not fluctuation around a line---to the contrary it
suggests two clear elbows in the CDF, which correspond to the two peaks.

%

%


%

It is also interesting to note that the empirical Fr\'echet mean was
observed to be insensitive to the choice of the bandwidth parameter
used in the construction of the estimated conditional mean measures
$\widehat\Lambda_i$. Of course, the warp functions $\widehat{T}_i$
themselves (and hence the registered processes) would depend on this
parameter, since these couple $\widehat\Lambda_i$ and $\widehat
{\lambda
}$---and while the latter is insensitive to the choice of bandwidth
parameter, the former is clearly not.

Further simulations carried out in the supplementary material \cite
{supplement} reaffirm these findings for different ``sample sizes''
$\tau
$ and choices of smoothing parameter. Furthermore, numerical
comparisons also carried out in the supplement suggest that Fr\'
echet--Wasserstein registration outperforms Fisher--Rao registration
(carried out as in \cite{anuj_arxiv} at the level of CDFs), in terms of
how close the registered processes are to the original point processes
(prior to warping), where ``closeness'' is measured by means of the
$\ell
_2$ distance of the ordered points. This is not surprising given our
unbiasedness considerations (Proposition~\ref{prop:umvue}), since the
Fisher--Rao estimator is generally not $d$-unbiased.

%
%
%

%

\subsection{Triangular Cox processes}
\begin{figure}

\includegraphics{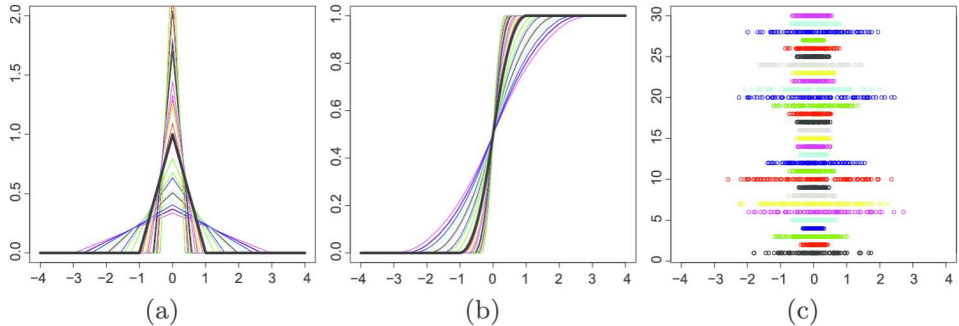}

\caption{\textup{(a)} Thirty triangular densities $f_{h_i}(t)$, with
$f_1$ in
solid black; \textup{(b)} Their corresponding distribution functions
$F_{h_i}(t)$, with $F_1$ in solid black; \textup{(c)} Thirty Cox processes,
constructed as follows: first, we generate $\Pi_i$ as i.i.d. Poisson
processes with mean density $f_1$, then we warp them by forming ${T_{i}}{}_{\#}
\Pi_i$.}
\label{fig:lineardendist}
\end{figure}

We now treat a second scenario that somewhat deviates from our model
assumptions, because it involves linear warp functions. Consequently,
phase variation can also be seen at the level of densities (see
Section~\ref{sec:measures_vs_densities}). Consider the family of
triangular densities of support length $2h$ and height $1/h$, and their
corresponding distribution functions (see Figure~\ref{fig:lineardendist})
\begin{eqnarray*}
f_h(t)&=&\frac{1}h \biggl(1-\frac{1}h\vert t \vert
\biggr),\qquad \vert t \vert\le h, h>0,\\
 F_h(t)&=& %
\cases{
\displaystyle\frac{1}{2h^2} (t+h )^2, &\quad $-h\le t\le0,$\vspace*{2pt}
\cr
\displaystyle 1-
\frac{1}{2h^2} (h-t )^2, &\quad$ 0\le t\le h$.} %
\end{eqnarray*}

Our example will consist in phase varying Poisson processes, with
structural mean distribution equal to $F_1$ (i.e., the triangular
distribution function with $h=1$). To this aim, let $h$ be a random
variable valued in $(0,C]$, so that the random measures have a common
support $I=[-C,C]$, but they are not strictly positive there. Following
the same steps as in the proof of Proposition~\ref{prop:uniqueness}, it
can be seen that the random measure with distribution function $F_h$
has a unique theoretical Fr\'{e}chet mean with distribution function
$F_{\mathbb E[h]}$, in the sense that for all distribution functions
$G\ne F_{\mathbb E[h]}$, we have (allowing for a slight abuse of
notation) $\mathbb{E} [d^2 (F_{\mathbb{E}[h]},F_{h} )
]<\mathbb{E} [d^2 (G,F_{h} ) ]$ (note
that Proposition~\ref{prop:uniqueness} and its proof remain valid as
long as the measures have no atoms; they do not need to be strictly
increasing). The warp map corresponding to an $h$ is $W_h(x)=hx$, and
it is not a homeomorphism of $I$ (unless $h=1$), thus violating our
assumptions (see Section~\ref{sec:measures_vs_densities}). To construct
our phase-varying point processes, we generate 30 i.i.d. copies $\{h_j\}
_{j=1}^{30}$ of a random variable\vspace*{1pt} $h$ following the mixture of uniform
distributions $\alpha\mathcal U[0.35,1]+(1-\alpha)\mathcal U[0.35,3]$,
where $\alpha=0.675$ is chosen so that $\mathbb E[h]=1$. Then we
generate 30 Poisson processes, with cumulative mean measure $\tau
\times
F_1$ [i.e., $h=1$, see Figure~\ref{fig:lineardendist}(c)], $\tau=93$, and
warp them by the maps $\{T_i=W_{h_i}\}_{i=1}^{30}$. This yields 30 Cox
processes, each with a realised directing measure {$93\times\Lambda
_1,\ldots,93\times\Lambda_{30}$, respectively, where the $\Lambda
_1,\ldots
,\Lambda_{30}$ }have distribution functions $F_{h_1},\ldots,F_{h_{30}}$
[depicted in Figure~\ref{fig:lineardendist}(b)]. The resulting warped
spike trains are displayed in Figure~\ref{fig:lineardendist}(c).

%

%

%

%

\begin{figure}

\includegraphics{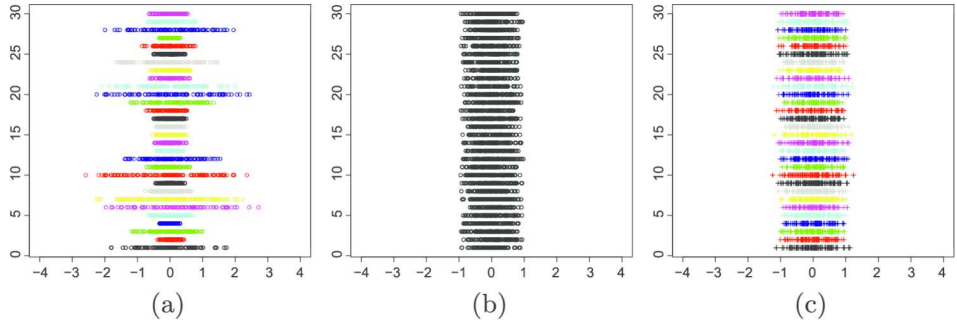}

\caption{Triangular Cox processes: \textup{(a)} The warped point
processes; \textup{(b)}
The original point processes; \textup{(c)} The registered point processes.}
\label{registration_plots_triangular}
\end{figure}

\begin{figure}

\includegraphics{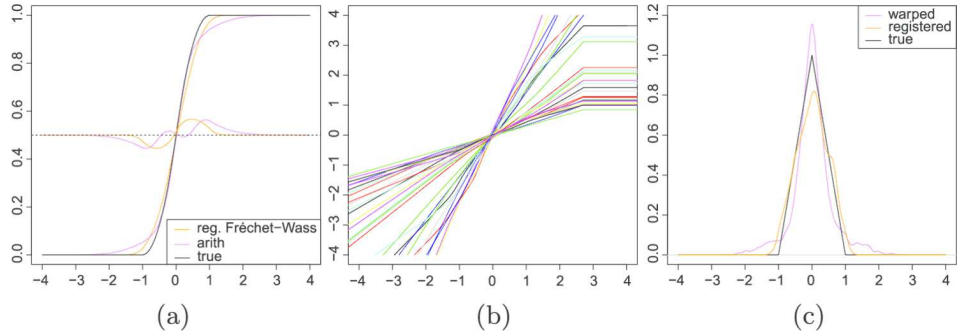}

\caption{\textup{(a)} The empirical arithmetic mean, our estimated regularised
Fr\'echet--Wasserstein mean, and the true mean CDF; \textup{(b)} The estimated
warp functions; \textup{(c)} Kernel estimates of the density function
of the
true structural mean, based on the original spike trains, and on the
registered spike trains.}\label{fig8}
\end{figure}

\begin{figure}[b]

\includegraphics{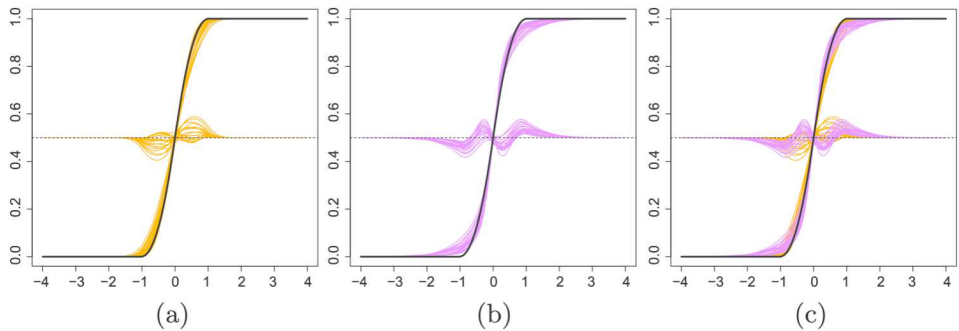}

\caption{\textup{(a)} Comparison of our estimated regularised Fr\'
echet--Wasserstein mean, and the true mean CDF, for 20 independent
replications of the experiment; \textup{(b)} Comparison of the
arithmetic mean,
and the true mean CDF, for the same 20 replications; \textup{(c)} Superposition
of \textup{(a)} and \textup{(b)}. In all three cases, the curves
oscillating about the
horizontal line $y=1/2$ are residual curves, centred at $1/2$.}
\label{triangle_replications}
\end{figure}

Assuming that the parametric form of the model is unknown to us, we
carry out the separation of amplitude and phase variation
nonparametrically, as described in Section~\ref{sec:estimation}. We
smooth each spike train using a Gaussian kernel with bandwidth chosen
by unbiased cross validation to obtain the estimators $\{\widehat
\Lambda
_i\}_{i=1}^{30}$ (strictly speaking, not in line with our discussion in
Section~\ref{sec:LambdaEstimation}, but this has no practical effect),
estimate the warp functions $\{\widehat T_i\}_{i=1}^{30}$, as described
in Section~\ref{sec:warp_estimation}, and produce a registration of the
point processes using these (Figure~\ref
{registration_plots_triangular}). We see that these warp functions
[Figure~\ref{fig8}(b)] are indeed nearly linear (besides
numerical instabilities at the boundary of the domain). The regularised
Fr\'{e}chet--Wasserstein mean of $\{F_{\widehat\Lambda_i}\}_{i=1}^{30}$
is depicted in Figure~\ref{fig8}(a), contrasted with the
arithmetic mean and the true structural mean. Note that the regularised
Fr\'echet--Wasserstein mean is supported on a subset of the domain, as
is the true structural mean; by contrast, the arithmetic mean is
supported almost on the entire domain, which is visible in Figure~\ref{fig8}(a), where it has left-and-right tails that persist.
Though both the regularised Fr\'echet--Wasserstein and the arithmetic
mean perform well near the point of symmetry of the structural mean
(which is to be expected, at least for the arithmetic mean, since the
location of the structural measure is invariant to the warp action),
the regularised Fr\'echet--Wasserstein mean estimates the support and
tails of the structural measure visibly better. These observations are
more clearly depicted in the residual plots contained in Figure~\ref
{triangle_replications}, where the residual curves of the deviation
between the arithmetic/Fr\'echet means and the estimand are considered,
for 20 independent repetitions of the same simulation experiment. It is
seen in that diagram that the arithmetic mean is clearly biased,
especially near the boundaries of the support of the true structural mean.

%

%

%

To gauge the effectiveness of the registration carried out, we also
constructed kernel estimators of the density of the structural mean,
based on the original (warped) point processes, and on the registered
(aligned) point processes. These are shown in Figure~\ref{fig8}(c). They illustrate that the density estimate based on
the raw data overestimates the mode as well as the tails of the true
density, whereas the density estimate based on the registered data fits
both the bulk and the tails of the~density quite nicely. As in the
previous example, the minor fluctuations of these density estimates are
not intrinsic to our registration procedure, but to the kernel density
estimation procedure.

Stability of the estimated structural mean CDF with respect to the
smoothing parameter was also observed in this example, and persisted in
additional simulations (presented in the supplementary material \cite
{supplement}), where different sample sizes were also considered.
Simulation comparisons showed that also in this scenario our approach
performs at least as well as the Fisher--Rao approach in terms of
registration of the point processes.

\section{Discussion}

We have introduced a framework formalising the confounding of amplitude
and phase variation in point process data, and demonstrated how this
can be used for their consistent nonparametric separation on the basis
of independent realisations thereof. The key ingredient of our approach
was the observation that for the point process warping problem, the
classical functional data assumptions on warp functions are equivalent
to the geometry of the Monge problem of optimal transportation.

A particularly attractive aspect of the present framework is that it
yields an identifiable setup, with a clear notion of over/under
registration through the concept of bias. Indeed, we prove that
consistent estimation of the warp functions \emph{is} possible in our
framework for point process data, circumventing the so-called
``pinching effect'' (see, e.g., Kneip and Ramsay \cite{kneip1},
Section~2.4) even under very sparse sampling regimes (Remark~\ref
{sparse_remark}). Furthermore, our consistency results present some
appealing features: there is no finite-dimensional parameterisation,
and the unknown warp functions and measures are allowed to be genuinely
functional, that is, infinite dimensional (contrary to, say Tang and
M\"uller \cite{tang}; Gervini and Gasser~\cite{gervini+gasser}; R\o nn
\cite{ronn}); though the consistency of the warp functions is in the
uniform metric, there is no need for the introduction of additional
smoothness penalties on the warp functions, and no tuning parameter
needs be selected to impose this (the regularity is inherited directly
from the underlying regularity of the structural and conditional mean
point process measures themselves; in the functional case, this
corresponds to the regularity of the curves themselves); consistency is
established with reference to a population, that is, the number of
``individuals'' (processes) is allowed to grow along with the ``density
of their sampling'' (with a clearly identified relationship between the
two), instead of establishing consistency conditional on the sample
(i.e.,  with a fixed number of curves, assuming only that the density of
sampling for each curve increasing, with no reference to a more general
``curve population,'' as in, e.g., Kneip and Engel \cite{kneipengel},
Wang and Gasser \cite{wang2}, and Gervini and Gasser \cite{selfmodel}).
In our experience, when consistency results are given in the functional
warping literature, they typically feature at least one of these
restrictions. We do not mention these characteristics as a claim to
superiority, but rather point them out as a special feature of the
problem in the point process case, afforded by the optimal
transportation geometry (since the very warping process is inextricably
linked with the metric structure of the space). Nevertheless, it is
interesting to note that the functional form of the warp function
estimator (\ref{warping_estimator}) is strikingly similar with the
pairwise synchronisation estimator of Tang and M\"uller \cite{tang}, equation
(7).

Further to consistency, we are able to obtain detailed rates of
convergence. These show $\sqrt{n}$-consistency and a central limit
theorem in the special case of warped Poisson processes (Cox processes)
under dense sampling. These can serve as a basis for uncertainty
quantification, but also indicate that our estimator can attain the
optimal rate of convergence under dense sampling.

Though we have demonstrated that the optimal transportation geometry is
canonical if warping occurs at the level of the spike train
observations (at the level of measures), it is possible to introduce
warping at the level of the density of the underlying mean measure (see
Section~\ref{sec:measures_vs_densities}). In such a framework, there
are options other than the optimal transportation geometry that be may
better suited for the formalisation of the warping problem. For
example, in the case of functional data, Tucker, Wu and Srivastava \cite
{fisher-rao} attack the warping problem by imbedding the data in a
quotient space modulo warp functions. This is done by employing a
Fisher--Rao-type metric, which is invariant with respect to the action
of a warping group. Recent work by Wu and Srivastava \cite{EJS} extends
their approach to the case of spike trains, by smoothing the spike
trains and considering them as densities in the Fisher--Rao space. This
geometry may be more natural than the optimal transportation one to
model phase variation at the level of densities.

A natural question for further work is that of \emph{multivariate phase
variation}. For example, is the ``canonicity'' of the optimal
transportation framework preserved, and can one fruitfully proceed in a
similar manner? The key challenge in this case is that, in the case of
measures on subsets of $\mathbb{R}^d$, $d>1$, evaluation of the
empirical Fr\'echet mean in closed form is impossible (see, e.g., Agueh
and Carlier \cite{agueh}). Approximations can be sought, for example,
via Gaussian assumptions (Cuturi and Doucet \cite{doucet}) or via
reduction to several 1D problems (Bonneel et al. \cite{boneel}).
Indeed, during the final preparation of this manuscript, we became
aware of interesting independent work in parallel by Boissard, LeGuic and Loubes
\cite{loubes}, who consider the problem of estimating Wasserstein
barycentres for measures on $\mathbb{R}^d$, and define ``admissible''
groups of deformations that mimic the 1D case, thus allowing for
consistent estimation and evaluation of the sample barycentre by
calculating successive means between pairs (i.e., by an iterated barycentre).

Finally, it should be mentioned that once phase and amplitude variation
have been separated, they could each be subjected to a further analysis
of their own. The amplitude variation clearly would be analysed by
means of \emph{linear PCA} tools, along the lines described in
Section~\ref{sec:amp_pp}. On the other hand, the phase variation can be
analysed by making further use of the geometrical properties described
in Section~\ref{sec:geometry}: for instance, via tangent space PCA
(see, e.g., Boissard, LeGuic and Loubes \cite{loubes}) or via geodesic PCA (see,
e.g., Bigot et al. \cite{bigot_pca}). Indeed, the form of the limiting
covariance function in our central limit theorem (Theorem~\ref
{thm:asynorm}) suggests that strong connections can be established
between Wasserstein PCA methodology and the separation of amplitude and
phase variation.

\section*{Acknowledgments}
This paper grew out of work presented at the Mathematical Biosciences
Institute (Ohio State University), during the\break \href
{http://mbi.osu.edu/event/?id=162}{``Statistics of Time Warping and
Phase Variation''} Workshop, November\break 2012. We wish to acknowledge the
stimulating environment offered by the Institute. We are grateful to an
Associate Editor and three referees for their insightful and
constructive comments. The paper has genuinely improved as a result of
the review process.

\begin{supplement}[id=suppA]
\stitle{``Amplitude and phase variation of point processes''\\}
\slink[doi]{10.1214/15-AOS1387SUPP} 
\sdatatype{.pdf}
\sfilename{aos1387\_supp.pdf}
\sdescription{The online supplement
contains more detailed simulation experiments.}
\end{supplement}

%
%

%



\printaddresses
\end{document}